\documentclass[12pt]{amsart}

\usepackage{amssymb,amscd,amsthm,verbatim,fancyhdr}
\usepackage{mathrsfs}

\begin{document}
\numberwithin{equation}{section}
\newtheorem{thm}{Theorem}
\newtheorem{lemma}{Lemma}[section]
\newtheorem{clm}[lemma]{Claim}
\newtheorem{remark}[lemma]{Remark}
\newtheorem{definition}[lemma]{Definition}
\newtheorem{cor}[lemma]{Corollary}
\newtheorem{prop}[lemma]{Proposition}
\newtheorem{statement}[lemma]{Statement}

%%  REMARK:  Writing [thm] on the lower lines makes all theorems, lemmas, etc. numbered in the same system as "Theorem".
%%  Take it out to have them each numbered separately.

\newcommand{\hdt}{{\dot{\mathrm{H}}^{1/2}}}
\newcommand{\hdtr}{{\dot{\mathrm{H}}^{1/2}(\mathbb{R}^3)}}
\newcommand{\R}{\mathbb{R}}
\newcommand{\ei}{\mathrm{e}^{it\Delta}}
\newcommand{\ltrt}{{L^3(\mathbb{R}^3)}}
\newcommand{\ldrd}{{L^d(\mathbb{R}^d)}}
\newcommand{\lprd}{{L^p(\mathbb{R}^d)}}
\newcommand{\lt}{{L^3}}
\newcommand{\ld}{{L^d}}
\newcommand{\lp}{{L^p}}
\newcommand{\rt}{\mathbb{R}^3}
\newcommand{\rd}{\mathbb{R}^d}
\newcommand{\X}{\mathfrak{X}}
\newcommand{\F}{\mathfrak{F}}
\newcommand{\hdhalf}{{\dot H^\frac{1}{2}}}
\newcommand{\hdthalf}{{\dot H^\frac{3}{2}}}
\newcommand{\hdo}{\dot H^1}
\newcommand{\rthmiz}{\R^3\times(-\infty,0)}
\newcommand{\q}[2]{{#1}_{#2}}
\renewcommand{\t}{\theta}
\newcommand{\lxt}[2]{L_{x,\,t}^{#1}}
\newcommand{\rr}{\sqrt{x_1^2+x_2^2}}
\newcommand{\ve}{\varepsilon}
\newcommand{\hdhrt}{\dot H^\frac{1}{2}(\mathbb{R}^3)}
\renewcommand{\P}{\mathbb{P} }
\newcommand{\RR}{\mathcal{R} }
\newcommand{\TT}{\overline{T} }
\newcommand{\e}{\epsilon }
\newcommand{\D}{\Delta }
\renewcommand{\d}{\delta }
\renewcommand{\l}{\lambda }
\newcommand{\To}{\TT_1 }
\newcommand{\ukt}{u^{(KT)} }
\newcommand{\ttil}{\tilde{T_1} }
\newcommand{\etl}{e^{t\Delta} }
\newcommand{\et}{\mathscr{E}_T }
\newcommand{\se}{\mathscr{E}}
\newcommand{\ft}{\mathscr{F}_T }
\newcommand{\eti}{\mathscr{E}^{\infty}_T }
\newcommand{\fti}{\mathscr{F}^{\infty}_T }
\newcommand{\xjn}{x_{j,n}}
\newcommand{\xjpn}{x_{j',n}}
\newcommand{\ljn}{\l_{j,n} }
\newcommand{\ljpn}{\l_{j',n} }
\newcommand{\lkn}{\l_{k,n} }
\newcommand{\voj}{U_{0,j} }
\newcommand{\uoj}{U_{0,j} }
\newcommand{\voo}{U_{0,1} }
\newcommand{\uon}{u_{0,n} }
\newcommand{\N}{\mathbb{N} }
\newcommand{\Z}{\mathbb{Z} }
\newcommand{\E}{\mathscr{E} }
\newcommand{\tu}{\tilde{u} }
\newcommand{\tU}{\tilde{U} }
\newcommand{\etj}{{\E_{T^*_j}}}
\newcommand{\B}{\mathcal{B}}
\newcommand{\tujn}{\tilde{U}_{j,n}}
\newcommand{\tujpn}{\tilde{U}_{j',n}}
\newcommand{\soj}{\sum_{j=1}^J}
\newcommand{\soi}{\sum_{j=1}^\infty}
\newcommand{\hnj}{H_{n,J}}
\newcommand{\enj}{e_{n,J}}
\newcommand{\pnj}{p_{n,J}}
\newcommand{\lfoi}{L^5_{(0,\infty)}}
\newcommand{\lfhoi}{L^{{5/2}}_{(0,\infty)}}
\newcommand{\wnj}{w_n^J}
\newcommand{\rnj}{r_n^J}
\newcommand{\lfij}{{L^5_{I_j}}}
\newcommand{\lfiz}{{L^5_{I_0}}}
\newcommand{\lfio}{{L^5_{I_1}}}
\newcommand{\lfit}{{L^5_{I_2}}}
\newcommand{\lfhij}{{L^{{5/2}}_{I_j}}}
\newcommand{\lfhiz}{{L^{{5/2}}_{I_0}}}
\newcommand{\lfhio}{{L^{{5/2}}_{I_1}}}
\newcommand{\lfhit}{{L^{{5/2}}_{I_2}}}
\newcommand{\lfi}{{L^5_{I}}}
\newcommand{\lfhi}{{L^{{5/2}}_{I}}}
\newcommand{\doh}{D^\frac{1}{2}}
\newcommand{\tjn}{t_{j,n}}
\newcommand{\tjpn}{t_{j',n}}
\newcommand{\wnlj}{w_n^{l,J}}
\renewcommand{\O}{\mathcal{O}}
\newcommand{\bes}{{\dot B^{s_p}_{p,q}}}
\newcommand{\bespp}{{\dot B^{s_p}_{p,p}}}
\newcommand{\besa}{{\dot B^{s_{a}}_{a,b}}}
\newcommand{\besb}{{\dot B^{s_{p}}_{p,q}}}
\newcommand{\besr}{{\dot B^{s_{b}+\frac{2}{\rho}}_{b,r}}}
\newcommand{\tn}{{\tilde \|}}
\newcommand{\tp}{{\tilde{\phi}}}
\newcommand{\xto}{\xrightarrow[n\to\infty]{}}
\newcommand{\LL}{\Lambda}
\newcommand{\binf}{{\dot B^{-{d/p}}_{\infty,\infty}}}
\newcommand{\brq}{{\dot B^{s_{p,r}}_{r,q}}}
\newcommand{\bfin}{{\dot B^{s_{p,r}}_{r,q}}}
\newcommand{\trnj}{{ \tilde{r}_n^J}}
\def\longformule#1#2{
\displaylines{ \qquad{#1} \hfill\cr \hfill {#2} \qquad\cr } }
\def\sumetage#1#2{ \sum_{\binom{\scriptstyle {#1}}{\scriptstyle {#2}}}}
\def\longformule#1#2{ \displaylines{\qquad{#1} \hfill\cr \hfill {#2} \qquad\cr } }

\title[Profile decomposition and regularity for N-S]{A profile decomposition approach to the $L^\infty_t(L^{3}_x)$ Navier-Stokes regularity criterion }
\author{Isabelle Gallagher}
\address{Institut de Math\'ematiques de Jussieu UMR 7586\\
Universit\'e Paris VII\\
175 rue du Chevaleret\\
75013 Paris, France}
\email{gallagher@math.univ-paris-diderot.fr}

\author{Gabriel S. Koch}
\address{
Oxford Centre for Nonlinear PDE\\
Mathematical Institute\\
University of Oxford\\
24-29 St Giles'\\
Oxford\\
OX1 3LB England 
}
\email{koch@maths.ox.ac.uk}
\author{Fabrice Planchon}
\address{Laboratoire J. A. Dieudonn\'e, UMR 6621\\
Universit\'e de Nice Sophia-Antipolis\\
Parc Valrose\\
06108 Nice Cedex 02\\
FRANCE }
\email{fabrice.planchon@unice.fr}
  \thanks{The first  author was partially supported by
the A.N.R grant ANR-08-BLAN-0301-01 "Mathoc\'ean", as well as  the Institut Universitaire de France.
The second author was supported by the EPSRC Science and Innovation award to the Oxford Centre for Nonlinear PDE (EP/E035027/1).
The third author was partially supported by A.N.R. grant SWAP}
\date{}
\maketitle
\begin{abstract}
In this paper we continue to develop an alternative viewpoint on
recent studies of Navier-Stokes regularity in critical spaces, a
program which was started in the recent work \cite{kk} by C. Kenig and
the second author.  Specifically, we prove that
 strong solutions which remain bounded in the space ${L^3(\mathbb{R}^3)}$ do not become
singular in finite time, a known result established in \cite{ess} by Escauriaza, Seregin and \v Sver\'
ak in the context of suitable weak solutions. Here, we use the method of ``critical elements" which was recently
developed by Kenig and Merle to treat critical dispersive equations. Our
main tool is a ``profile decomposition" for the Navier-Stokes
equations in critical Besov spaces which we develop here.  As a
byproduct of this tool, assuming a singularity-producing initial datum
for Navier-Stokes exists in a critical Lebesgue or Besov space, we
 show there is one with minimal norm, generalizing a result of Rusin
and \v Sver\' ak~\cite{sverakrusin}.
\end{abstract}
\section*{Introduction}
We consider the incompressible Navier-Stokes equations in $\mathbb{R}^d$,
\begin{equation}
\label{ns}
(NS) \quad \left \{ \begin{array}{rcl} \displaystyle{\frac{\partial u}{\partial
t}} & = & \Delta u - \nabla \cdot (u\otimes u)-\nabla \pi,\\ \nabla
\cdot u & = & 0,\\ u|_{t=0} & = & u_0
\end{array} \right.
\end{equation}
for $(x,t)\in \rd \times (0,T)$, where $u=u(x,t)$ is the velocity vector field and~$\pi(x,t)$ is the associated pressure function.
For $d\geq 3$, global weak solutions are known to exist, but their
uniqueness (as well as the smoothness of the solution for smooth data) has remained an open problem since the pioneering work
\cite{leray}. There exist several conditional results, of which
Serrin's criterion  is perhaps the most well-known: if a weak solution $u$
is such that
\begin{equation}
  \label{eq:serrin}
  u\in L^p([0,T];L^q(\mathbb{R}^d)) \text{ with } \frac 2 p +\frac d q=1,\,\,q>d,
\end{equation}
then $u$ is smooth on $(0,T)$. On the other hand, there is a long line of work on
constructing local in time solutions, from \cite{KF} to \cite{KT}. In
this framework of local in time (strong, e.g. unique) solutions,
Serrin's criterion may be understood as a non blow-up criterion at time $T$: e.g. if
$u$ is a strong solution with $u_0\in L^d(\mathbb{R}^d)$, that is
$u\in C([0,T);L^d(\mathbb{R^d}))$, and if \eqref{eq:serrin} is satisfied,
then one may (continuously and uniquely) extend the solution $u$ past time $T$.

In the recent important work \cite{ess}, Escauriaza-Seregin-\v Sver\'ak
obtained the endpoint version of Serrin's criterion\footnote{\cite{ess} treats the case $d=3$; the case $d>3$ was treated similarly later in \cite{dongdu2}.}: $u\in
L^\infty([0,T];L^3(\mathbb{R}^3))$ implies no blow-up; they work with
the so-called suitable weak solutions introduced in \cite{CKN}. A similar, but much more precise result  in the (smaller) space~$
L^\infty([0,T];\dot H^\frac12 (\mathbb{R}^3))$ was obtained recently by Seregin in~\cite{seregin}; the main point in that result is that it is proved that  the $\dot H^\frac12$ norm of the solution $u(t)$ blows up  as~$t$ goes to blow-up time,  and not just for a subsequence as is known in the~$L^3(\mathbb{R}^3) $ case (see~\cite{seregin0} for a partial result in that direction)\footnote{The result~\cite{seregin0} was actually very recently improved by Seregin in~\cite{seregin1}, who obtains the blow up of~$u(t)$ in the~$L^3(\mathbb{R}^3)$ case.}.

Our goal here is to obtain such a result by a somewhat different route,
following the concentration-compactness methods developed by
Kenig-Merle in the context of energy critical dispersive equations (\cite{KMNLS,
  KMNLW}) and then extended to subcritical problems~(\cite{KMNLS3}) or
supercritical problems~(\cite{KMNLW7}). In our context, the
Navier-Stokes equations are supercritical with respect to their only
known a priori bound, which is the energy inequality (the
$L^2(\mathbb{R}^3)$ norm of a solution is decreasing), while a scale
invariant norm for the data is $\dot H^\frac 1 2(\mathbb{R}^3)$ or
$L^3(\mathbb{R}^3)$. In \cite{kk}, C. Kenig and the second author
carried out such a program for solutions~$u\in L^\infty([0,T];\dot
H^\frac 1 2(\mathbb{R}^3))$. The first step in following the
Kenig-Merle roadmap is to prove the existence of a so-called
``critical element'' (or minimal blow-up solution), which follows from
suitable profile decompositions of the data, the linear solution (to the heat equation) and
the nonlinear solution (to the full Navier-Stokes system). Such decompositions were introduced by P. G\'erard in
\cite{PG} to study the defect of compactness of the Sobolev embedding
theorem, and then used by the same author and H. Bahouri to study the critical defocusing wave equation
\cite{BG}. In the context of Navier-Stokes, they were developed in~\cite{Gallagher} by the first author, and served as a crucial tool in
implementing the roadmap in \cite{kk}; this explains why the result
in \cite{kk} applies to~$\dot H^\frac 1 2(\mathbb{R}^3)$ rather than
$L^3(\mathbb{R}^3)$, as \cite{Gallagher} implements profile decompositions
in the Sobolev scale. Recently the second author extended profile
decompositions (\cite{gk}) to study the embedding  $ L^{d}(\mathbb{R}^d) \hookrightarrow
\dot{B}^{\frac{d}{p}-1}_{p,q}(\mathbb{R}^d)$, where the latter are
((NS)-critical) homogeneous Besov spaces, with $p>d$, or more generally within the
Besov scale itself.

Our main goal is threefold:
\begin{itemize}
\item  we develop  profile decompositions for solutions to the
  Navier-Stokes equations with data in  $ L^{d}(\mathbb{R}^d)$ or
  $\dot{B}^{\frac{d}{p}-1}_{p,q}(\mathbb{R}^d)$ with~$1\leq p,q<2d+3$,
  extending the results from \cite{Gallagher}; the main difficulty here compared to \cite{Gallagher} consists in handling multilinear
  interactions between profiles and remainders,  as well as a lack of
  orthogonality in $L^3(\mathbb{R}^3)$ for the profile decomposition.
\item we use this profile decomposition to implement the
  Kenig-Merle roadmap for solutions $u\in
  L^\infty([0,T];L^3(\mathbb{R}^3))$. The
  Hilbert nature of $\dot H^\frac 1 2(\mathbb{R}^3)$ proved helpful  in \cite{kk} to deal    with weak convergence issues,
  as well as again
 with  multilinear interactions.  We need to face these issues here, however eventually  we obtain a streamlined argument which
  leads to weak convergence toward zero at blow-up time for critical
  elements. Once this ``compactness'' result is proved, Serrin's
  endpoint criterion is obtained as in \cite{kk}, following closely
  the backward uniqueness argument of~\cite{ess};
\item  we use the profile decomposition in another direction,
  extending recent work of Rusin-\v Sver\'ak \cite{sverakrusin}: we prove that there
  always exists a minimal blow-up initial datum in $L^d(\mathbb{R}^d)$
  if any such datum exists, and that moreover the set of such data is
  compact in~$\dot{B}^{\frac{d}{p}-1}_{p,q}(\mathbb{R}^d)$, with $d<
  p\leq q\leq +\infty$, up to transformational invariance of the
  equations. We moreover prove a similar statement, involving two
  different Besov spaces in the aforementioned scale, with $p<2d+3$.
\end{itemize}
The next section introduces the function spaces we shall be using, and collects a few  well-known facts
about the Navier-Stokes system and its solutions in those function spaces.  The next three sections are devoted to the
profile decomposition, the regularity criterion and the minimal
blow-up data, respectively.\\

{\bf Acknowledgement:}  The second author would like to express his sincere thanks to Professor Carlos Kenig for suggesting to him the problem which we treat in Section 3 below.\\

After completion of this work, we learned of~\cite{jiasverak} where  a result in  the same spirit as our last section is proved, namely the existence of initial data with minimal~$L^3(\R^3)$ norm for potential Navier-Stokes singularities; in~\cite{jiasverak}  the compactness in~$L^3(\R^3)$ up to translation-dilation is also obtained.

\section{Preliminaries}

For the convenience of the reader, we start by recalling the usual definition of
Besov spaces. We usually write $X$ as a shorthand for
the function space $X(\mathbb{R}^d)$, where $d$ is the space dimension.
\begin{definition}
\label{d1}
Let $\phi$ be a function in $\mathcal{S} $ such that $\widehat\phi =
1$ for $|\xi|\leq 1$ and $\widehat\phi= 0$ for $|\xi|>2$, and define
$
\phi_{j}(x):= 2^{dj}\phi(2^{j}x).$
Then the frequency localization operators are defined by
$$
S_{j} := \phi_{j}\ast\cdot, \quad
\Delta_{j} := S_{j+1} - S_{j}.
$$
Let $f$ be in $\mathcal{S}' $. We say $f$ belongs to~$\dot
B^{s}_{p,q}$ if
\begin{itemize}
\item The partial sum $ \sum^{m}_{-m} \Delta_{j} f $ converges to $f$
  as a tempered distribution if~$s < {\frac{d}{p}}$ and after taking
  the quotient with polynomials if not.
\item The sequence $\epsilon_{j} := 2^{js}\| \Delta_{j}  f \|_{L^{p}}$
  belongs to $\ell^{q}$, and its $\ell^q$-norm defines the Besov norm of $f$.
\end{itemize}
\end{definition}
We shall also need a slight modification of those spaces, taking into
account the time variable.
\begin{definition}
  \label{raah11}
  Let $u(x,t)\in \mathcal{S}'(\mathbb{R}^{d+1})$ and let $\Delta_j$ be a
  frequency localization with respect to the $x$ variable. We shall say
  that $u$ belongs to~$ {\mathcal{L}^{\rho}}([ a,b ] ;\dot{B}^{s}_{p,q})$ if $u(t)\in \dot B^s_{p,q}$ for all $t\in [a,b]$ and
  \begin{equation*}
  %\label{eq:raaah12}
  2^{js}\|\Delta_j u\|_{L^\rho([ a,b ];L^{p}_x)} =:\varepsilon_j \in \ell^q\ .
\end{equation*}
We define
$$
\|u\|_{{\mathcal{L}^{\rho}}([ a,b ] ;\dot{B}^{s}_{p,q})} :=  \|2^{js}\|\Delta_j u\|_{L^\rho([ a,b ];L^{p}_x)} \|_{\ell^q},
$$
and $ {\mathcal{L}^{\rho}_T}(\dot{B}^{s}_{p,q}) :=  {\mathcal{L}^{\rho}}([ 0,T ] ;\dot{B}^{s}_{p,q})$.
\end{definition}
\begin{remark}
In the case where $ \rho \geq q$ one has of course the embedding~$ {\mathcal{L}^{\rho}}([ a,b ] ;\dot{B}^{s}_{p,q}) \hookrightarrow L^{\rho}( [ a,b ] ;\dot{B}^{s}_{p,q}) $ due to Minkowski's inequality.
\end{remark}
Let us introduce some notation. For any~$p$ in~$[1,\infty) $ we define
$$
s_p := -1 + \frac d p \cdotp
$$
For any initial datum $u_0\in \dot{B}^{s_p }_{p,q}$, with
$d<p\leq q<+\infty$, we shall denote by $NS(u_0)$ the local in time
strong solution to the Navier-Stokes equation \eqref{ns}.  For clarity,  by ``solution" to \eqref{ns} in the strong (sometimes called ``mild") sense, we mean a divergence-free solution $u$ to
$$
u_t = \D u - \P \nabla \cdot (u\otimes u)\ , \qquad u|_{t=0} = u_0
$$ (equivalent to solving \eqref{ns} for the ``right" $\pi$) in the Duhamel sense, where $\P$ is the projection operator onto divergence-free vector fields. Such
solutions were obtained in \cite{Cannone} for $3<p\leq 6$ and for all
$p<+\infty$ in \cite{Planchon}, and we refer to the appendix of
\cite{gip3} for a proof which is taylored to our purposes. The specific
case of $L^d (  \hookrightarrow \dot{B}^{-(1-\frac d p)}_{p,q}) $ data is included in such a result, as any additional
``regularity'' is propagated along the flow (see again \cite{gip3} for
a   proof of this well-known fact).

\medskip

Let us make those results on the Cauchy problem more precise. We define the function space
$$
E_{p,q} (T) := \mathcal{L}^\infty([0,T];{\dot B^{s_{p}}_{p,q}}) \cap
\mathcal{L}^{\frac{2p}{p+1}} ([0,T];\dot B^{s_{p}+1+\frac 1 p}_{p,q})\ .
$$
In particular, $E_{p,q} \subset \mathcal{L}^r(\dot B^{s_p + \frac{2}{r}}_{p,q})$ for $\frac{2p}{p+1} \leq r \leq \infty$ by interpolation.  We recall (see e.g.~\cite{gip3})  that~$NS(u_0)$ belongs to~$E_{p,q} (T) $ for some time~$T$, and one may define a maximal time
$T^* = T^*(u_0)$ such that this holds for any $T<T^*$ (actually the solution belongs to~$ \mathcal{L}^\infty([0,T];{\dot B^{s_{p}}_{p,q}}) \cap
\mathcal{L}^{1} ([0,T];\dot B^{s_{p}+2}_{p,q})$ but that fact will not be used here). If the initial
datum is small enough then~$T^* = \infty$ (and under such a condition one may include $q=\infty$, although one cannot in general obtain local solutions for $q=\infty$). Moreover, $u$ belongs to~$E_{p,q} (T^*)$ if and only if~$T^* = \infty$, and in that case one has (see~\cite{gip3}) that~$\displaystyle \lim_{t \to \infty} \|NS(u_0)(t)\|_{\dot
B^{s_{p}}_{p,q}} = 0$ as well.
Finally recall that if~$NS(u_0)$ belongs
to~$E_{p,q} (T) $ and if~$u_0$ belongs to~${\dot B^{s_{a}}_{a,b}}$
(resp. ${L^d(\mathbb{R}^d)}$) with~$a \leq p$ and~$b \leq q$,
then~$NS(u_0)$ belongs to~$E_{a,b} (T) $
(resp. $C([0,T];{L^d(\mathbb{R}^d)}$) with the same life span (see \cite{gip3} for instance, or~\cite{adt}).
  \section{Profile Decompositions}\label{profile}
 \subsection{Notation and statement of the result}
In what follows, we shall need the following notion, where the dimension~$d$ is always chosen such that~$d \geq 2$.
 \begin{definition}\label{deforthseq} For $j\in \N$, let $\{(\lambda _{j,n} ,x_{j,n})\}_{n=1}^\infty \subset (0,\infty) \times \mathbb{R}^d$ be a sequence of ``scales" $\ljn$ and ``cores" $\xjn$.  We say that such a set of sequences
is (pairwise)  orthogonal if
\begin{equation}\label{orthseq} \quad j \neq j' \quad \Longrightarrow
\quad \left\{ \quad
\begin{array}{c} \textrm{either} \quad \displaystyle{\lim_{n \to
+\infty} \frac{\lambda _{j,n} }{\lambda _{j',n} } + \frac{\lambda
_{j',n} }{\lambda _{j,n} } = +\infty}\\\\ \textrm{or} \\\\
\displaystyle{ \frac{\lambda _{j,n} }{\lambda _{j',n} } \equiv 1 \quad
\textrm{and} \quad \lim_{n \to +\infty} \frac{|x_{j,n} -
x_{j',n}|}{\lambda _{j,n} } = +\infty\ .}
\end{array} \right .
\end{equation}
\end{definition}
Let us define, for functions $f \in {L^d(\mathbb{R}^d)}$,
norm-invariant transformations depending on translation and scaling
parameters $x_{j,n}$ and $\lambda_{j,n}$ respectively by
\begin{equation}
  \label{eq:rescaling}
  \Lambda_{j,n} f(x):= \frac{1}{\lambda _{j,n} } f \left(\frac{x-x_{j,n}}{\lambda _{j,n} }\right)\ .
\end{equation}
Then the reason for the terminology in the previous definition becomes
clear by considering the following quantity:
\begin{multline}\label{Ldorthogx}
\int_{\mathbb{R}^d} |\Lambda_{j,n}f(x)|^{d-1} |\Lambda_{j',n}g(x)|\ dx
\\
= \frac{\lambda _{j,n}}{\lambda _{j',n}} \int_{\mathbb{R}^d}
|f(y)|^{d-1} \left|g\left(\frac{\lambda _{j,n}}{\lambda
_{j',n}}y + \frac{x_{j,n} - x_{j',n}}{\lambda _{j',n}}\right)\right|\
dy\ .
\end{multline}
One similarly has equality if on the right-hand side one interchanges~$j$ and~$j'$ and the arguments of the functions, and hence such a term tends to
zero as~$n\to\infty$ if (\ref{orthseq}) holds, since we may
approximate~$f$ and $g$ in~${L^d(\mathbb{R}^d)}$ by~$\mathcal{C}^\infty_0$-functions.  This implies that
$$\|\Lambda_{j,n}f+\Lambda_{j',n}g\|^d_{L^d(\mathbb{R}^d)} =
\|\Lambda_{j,n} f\|^d_{L^d(\mathbb{R}^d)} +\|\Lambda_{j',n} g\|^d_{L^d(\mathbb{R}^d)} + \circ(1)$$
as~$n\to\infty$ if the scales are orthogonal.  Therefore
$\Lambda_{j,n}f$ and $\Lambda_{j',n}g$ are ``asymptotically
orthogonal" in ${L^d(\mathbb{R}^d)}$ (and similar statements will be
shown to hold in other critical spaces as well).

We first recall two theorems which were proved\footnote{It
is known that the wavelet-basis characterization of scalar
function-spaces used in \cite{gk} extend as well to vector-fields (and
in fact one may use divergence-free wavelets, see e.g. \cite{bf}).
Moreover, the slightly more specific formulations we give here are a
simple consequence of the theorems in \cite{gk}.  } in \cite{gk}. The first one deals with the defect of compactness of the embedding~$L^d  \hookrightarrow \dot B^{s_p}_{p,q}$ with~$ d < p,q \leq +\infty$.
\begin{thm}[\cite{gk}]\label{thm:dataprofa}
Let $p,q
\in \mathbb{R}$ satisfy $ d < p,q \leq +\infty$.  Let~$\{\varphi_n\}_{n=1}^\infty$ be a bounded
sequence in~${L^d(\mathbb{R}^d)}$ and let $\phi_0$ be any weak limit
point of~$\{\varphi_n\}_n$.  Then, after possibly replacing~$\{ \varphi_n\}_n$ by a subsequence which we relabel~$\{\varphi_n\}_n$, there exists a sequence of profiles~$\{\phi_j\}_{j=1}^\infty$ of~${L^d(\mathbb{R}^d)}$, and for each~$j\in \mathbb{N} $ sequences~$\{(\lambda _{j,n}
,x_{j,n})\}_{n=1}^\infty $ in~$ (0,\infty) \times \mathbb{R}^d$ of
scales and cores which are orthogonal in the sense of Definition~\ref{deforthseq}
such that, for all~$n,J \in \mathbb{N}$, if we define~$\psi_n^J$ by
\begin{equation}\label{profiles} \varphi_n(x) = \phi_0(x) +
\sum_{j=1}^J \frac{1}{\lambda _{j,n}}
\phi_j\left(\frac{x-x_{j,n}}{\lambda _{j,n}}\right) + \psi_n^J(x)
\end{equation}
the following properties hold:
\begin{itemize}
\item the function $\psi_n^J$ is a remainder in the sense that
\begin{equation}\label{orth2} \lim_{J\to\infty} \left(\limsup_{n\to
\infty} \|\psi_n^J\|_{\dot B^{s_p}_{p,q}} \right) = 0\,;
\end{equation}
\item there is a norm\footnote{This norm is based on a wavelet basis
expansion for functions in Besov and Triebel-Lizorkin spaces.} $\|
\cdot \tilde \|_{L^d(\mathbb{R}^d)}$ which is equivalent to $\|\cdot
\|_{L^d(\mathbb{R}^d)}$ such that
\begin{equation}\label{orth3}
\sum_{j=1}^\infty \|\phi_j{\tilde
\|}_{L^d(\mathbb{R}^d)}^d \leq \liminf_{n\to\infty}
\|\varphi_{n}{\tilde \|}_{L^d(\mathbb{R}^d)}^d
\end{equation}
and, for each integer~$J$,
\begin{equation}\label{orth4} \|\psi_n^J {\tilde
\|}_{L^d(\mathbb{R}^d)} \leq \|\varphi_{n}{\tilde
\|}_{L^d(\mathbb{R}^d)} + \circ(1) \quad \textrm{as} \quad n\to\infty\,;
\end{equation}
\item for any integer~$j $,  the following properties hold: either $\lambda _{j,n} = 1$ and
${\displaystyle \lim_{n\to\infty} |x_{j,n}| = +\infty}$, or
${\displaystyle \lim_{n\to \infty} \lambda _{j,n} \in \{0,+\infty\}}$.
\end{itemize}
\end{thm}
The second statement   deals with the
 defect of compactness of the embedding~${\dot
B^{s_{a}}_{a,b}}  \hookrightarrow \dot B^{s_{p}}_{p,q} $ with~$1\leq a < p \leq +\infty$ and~$1\leq b \leq (p/a)b    \leq q \leq
+\infty$.
\begin{thm}[\cite{gk}]\label{thm:dataprofb}
 Let $a, b,p,q\in [1,\infty]$ satisfy $1\leq a < p \leq +\infty$ and~$1\leq b \leq (p/a)b    \leq q \leq
+\infty$. Let~$\{\varphi_n\}_{n=1}^\infty$ be a bounded sequence in~${\dot
B^{s_{a}}_{a,b}}(\mathbb{R}^d)$, and let~$\phi_{0}$ be any weak limit
point of~$\{\varphi_n\}$. Then, after possibly replacing~$\{ \varphi_n\}_n$ by a subsequence which we relabel~$\{\varphi_n\}_n$, there exists a sequence of profiles~$\{\phi_j\}_{j =1}^\infty $ in~${\dot B^{s_{a}}_{a,b}}$, and for
each integer~$j\geq 1$ a sequence $\displaystyle\{(\lambda _{j,n}
,x_{j,n})\}_{n=1}^\infty$ of scales and cores which are orthogonal in
the sense of Definition~\ref{orthseq} such that, for all~$n,J\in \mathbb{N}$, if we define $\psi_n^J$ by
\begin{equation}\label{profilesa} \varphi_n(x) = \phi_{0} (x)+
\sum_{j=1}^J \frac 1{\lambda _{j,n} } \phi_j\left(\frac{ x - x_{j,n}}{
\lambda _{j,n} } \right) + \psi_n ^J(x)
\end{equation}
the following properties hold:
\begin{itemize}
\item the function $\psi_n^J$ is a remainder in the sense that
\begin{equation}\label{orth2a} \lim_{J\to\infty}
\left(\limsup_{n\to \infty} \|\psi_n^J\|_{\dot B^{s_{p}}_{p,q}}
\right) = 0\,;
\end{equation}
\item there is a norm $\| \cdot \tilde \|_{\dot B^{s_{a}}_{a,b}}$
which is equivalent to $\|\cdot \|_{\dot B^{s_{a}}_{a,b}}$ such that
for each integer~$n \in \mathbb{N} $, denoting $\tau := \max \{a,b\}$,
one has
\begin{equation}\label{orth3a} \left\| \left(\| \phi_j {\tilde
\|}_{\dot B^{s_{a}}_{a,b}} \right)_{j=0}^\infty \right\|_{\ell^\tau}
\leq \liminf_{n'\to\infty} \|\varphi_{n'}{\tilde \|}_{\dot
B^{s_{a}}_{a,b}}
\end{equation} and, for each integer~$J $,
\begin{equation}\label{orth4a} \| \psi_n^J {\tilde \|}_{\dot
B^{s_{a}}_{a,b}} \leq \|\varphi_{n}{\tilde \|}_{\dot B^{s_{a}}_{a,b}}
+ \circ(1) \quad \textrm{as} \quad n\to\infty\,;
\end{equation}
\item  for any integer~$j $,  the following properties hold: either $\lambda _{j,n} = 1$ and
${\displaystyle \lim_{n\to\infty} |x_{j,n}| = +\infty}$, or
${\displaystyle \lim_{n\to \infty} \lambda _{j,n} \in \{0,+\infty\}}$.
\end{itemize}
\end{thm}
It should be noted (see e.g. the proof of (\ref{wkzero}) below with $s=0$) that the~$ \phi_j$'s
are weak limit points of the sequence
$$
\displaystyle \lambda _{j,n} \varphi_n \left( \lambda _{j,n} \cdot +
x_{j,n} \right) .
$$
Consequently if Theorem~\ref{thm:dataprofa} or Theorem
\ref{thm:dataprofb} is applied to a sequence of divergence free vector
fields, each profile produced by the theorem is also a divergence free
vector field.

We shall now apply these results to a sequence of
bounded initial data to~(\ref{ns}). The following statement is the analogue
of Theorem 2 of~\cite{Gallagher} in the context of critical Besov
spaces. For the sake of simplicity, we shall only consider the case
$p=q$ ; the restriction on $p$ is technical in nature but it allows to
deal with $L^3(\R^3)$ data, and we plan to address this issue in future work.
\begin{thm}[NSE Evolution of Profile Decompositions]\label{thm:profa}
Suppose $d\geq 3$ and fix $a,b,p \in \R$ satisfying $d < a \leq b < p < 2d+3$. Let~$\{\varphi_n\}_{n=1}^\infty$ be a bounded sequence of
divergence-free vector fields in~${\dot B^{s_{a}}_{a,b}}(\rd)$, and let~$\phi_{0}$ be any weak limit point of~$\{\varphi_{n}\}$. Then, after
possibly relabeling the sequence due to the extraction of a
subsequence following an application of Theorem~\ref{thm:dataprofb}
(of which we retain the same notation for the profile decomposition), defining $u_n:=
NS(\varphi_{n})$ and~$U_j := NS(\phi_j)\in E_{a,b}(T_j)$ for every
integer~$j \in \mathbb{N} $ (where~$T_j $ is any real number smaller
than the life span~$T^*_j$ of~$\phi_j$, and~$T_j = \infty$ if~$T^*_j =
\infty$), the following properties hold:
\begin{itemize}
\item there is a finite (possibly empty) subset~$I $ of
$\mathbb{N}  $ such that\footnote{This is equivalent to $T^*_j < \infty \iff j\in I$ by \cite{gip3} and the small-data theory.}
$$
\forall j \in I, \quad T_j < \infty \quad \mbox{and} \quad \forall j
\in \mathbb{N}  \setminus I , \quad U_j \in
E_{a,b}({\infty}).
$$
Moreover  setting $\displaystyle \tau_n := \min_{j \in I}\lambda _{j,n} ^2
T_j $ if $I$ is nonempty and $\displaystyle \tau_n := \infty$ otherwise, we have  $$\sup_n \|u_n\|_{E_{a,b} ({\tau_n})}<\infty.$$

\item if there exists a time~$T \in
\mathbb{R}^+ \cup \{+\infty\}$ such that~$\{u_n\}$ is  uniformly bounded
in~$E_{a,b}(T)$, then if $I$ is nonempty we must have
\begin{equation}
\label{eq:Tsmall} \forall n \in \mathbb{N} , \quad T < \min_{j \in
I} \: \lambda _{j,n} ^2 T^*_j,
\end{equation}
and therefore in such a case the scales of concentration for which  $\displaystyle \lim_{n \rightarrow
\infty} \lambda _{j,n} = 0$ (small scales) generate global solutions
of~$(NS)$ (i.e., the corresponding $T^*_j =+ \infty$).

\item finally there exists some large $J_0 \in \N$ such that for
each $J >J_0$, there exists $N(J)\in \N$ such that for all~$n >N(J)$, all~$t \leq \tau_n$ and all~$x \in \mathbb{R}^d$, setting~$ w_n^J := e^{t\Delta}(\psi_n^J)$ and defining~$r_n^J$ by
\begin{equation}\label{evolprof} u_n (x,t) = U_0(x,t) + \sum_{j =
1}^{J} \frac{1}{\lambda _{j,n} } U_j \left( \frac{x - x_{j,n}}{\lambda
_{j,n} }, \frac{t}{\lambda _{j,n} ^2}\right) + w_n^J (x,t) + r_n^J
(x,t),
\end{equation}
then $w_n^J$ and $r_n^J$ are small remainders in the sense that
\begin{equation}\label{remainders} \lim_{J\rightarrow \infty}
\left(\limsup_{n \rightarrow \infty}
\|w_n^J\|_{E_{p,p}(\infty)}\right) = \lim_{J\rightarrow \infty}
\left(\limsup_{n \rightarrow \infty}
\|r_n^J\|_{E_{p,p}(\tau_n)}\right) = 0\ ;
\end{equation}
\end{itemize}
 \end{thm}
 \begin{remark}
 As a consequence of the proof of Theorem~\ref{thm:profa}, the
 decomposition~\eqref{evolprof} actually holds for all indices~$J$
 (once the theorem is proved as stated, the remainder $r_n^J$ for
 $1\leq J \leq J_0$ may be defined by~\eqref{evolprof} which now holds up to
 time $\tau_n$) but it should be clear that such a
 decomposition is of interest mostly when enough profiles have been
 extracted, yielding a small remainder.

 \end{remark}

  \begin{remark}\label{rkaaa}
  Although the first step in proving
Theorem \ref{thm:profa} is to apply Theorem \ref{thm:dataprofb}, if
moreover $\{\varphi_n\}_{n=1}^\infty$ is a bounded sequence of
divergence-free vector fields in ${L^d(\mathbb{R}^d)}$, one may
instead first apply Theorem \ref{thm:dataprofa} and standard
embeddings to see that the sequence enjoys all the properties in the
conclusion of Theorem \ref{thm:dataprofa} as well as those in the
conclusion of Theorem \ref{thm:profa}.
\end{remark}

\begin{remark}\label{remarkindices}
The   bound $p<2d+3$ is directly related to
  Lemma~\ref{lemmadrift} below; we did not seek the optimal values
  (w.r.t. the original $a,b$), and for $a=b=d$, the regularity
  $s_p+2/p$ which appears in this Lemma may be taken positive by
  choosing $p<d+2$. This in turn would avoid the direct use of
  paraproduct estimates in the proof of Lemma \ref{lemmaremainder},
  but we feel worth pointing out that negative regularity can be
  allowed in Lemma \ref{lemmadrift}. At any rate, paraproduct
  techniques are at the heart of the estimates in \cite{gip3} or any
  of the product estimates, both of which we heavily rely on.
\end{remark}

\subsection{Proof of Theorem~\ref{thm:profa} } The first step consists
of course in appying Theorem~\ref{thm:dataprofb} (or Theorem \ref{thm:dataprofa}, if applicable) to the sequence of
initial data. We briefly comment on the choice of indices: our main
goal being to deal with a sequence of data in $L^3(\R^3)$, a natural
choice (in spatial dimension $d$) is~$a=b=d$, for which in the profile
decomposition of the data we may pick any $p=q$ close to $d$ (hence
even less than $d+2$, see Remark~\ref{remarkindices}). For general $a,b$, we may relabel
$a=b=\sup(a,b)$ (using Sobolev's embedding), and then use Theorem \ref{thm:dataprofb} for $p=q$
strictly larger than $a$. Thus we can always use
$p=q=\sup(a,b)+\varepsilon$ as our indices for the remainder space.  Now let us continue with the proof.

\medskip

  With the notation of Theorem~\ref{thm:dataprofb}  we define
 $$
 u_n:= NS(\varphi_{n}) , \quad U_j := NS(\phi_j)\in E_{a,b}(T_j)\quad
\mbox{and} \quad w_n^J := e^{t\Delta}(\psi_n^J) .
 $$
 Then due to~\eqref{orth2a} and standard linear heat estimates we have
 \begin{equation}
   \label{eq:heatw}
   \lim_{J\rightarrow \infty} \left(\limsup_{n \rightarrow \infty}
\|w_n^J\|_{E_{p,p}(\infty)}\right) = 0\ .
 \end{equation}
 Moreover due to the orthogonality property~\eqref{orth3a}, the
sequence~$\displaystyle \{ \phi_j\}$ goes to zero in the space~${\dot B^{s_{a}}_{a,b}}$ as~$j$ goes to
infinity.  This implies that
there is~$j_0 $ such that for all~$j > j_0$, there is a global unique
solution associated with~$\phi_j$, as $\|\phi_j\|_{\dot
  B^{s_a}_{a,b}}<\varepsilon_0$ (the smallness constant of small data
theory). Hence,~$I$ will be a subset of $\{0,\dots,j_0\}$ which
proves the first part of the first statement in Theorem~\ref{thm:profa}.

All other statements will follow
  from obtaining careful bounds on both
profiles and $r_n^J$, by taking advantage of the local Cauchy theory
and its perturbed variants (as set up e.g. in \cite{gip3}; see the appendix for a slightly more general statement).

By the local Cauchy theory we can solve the Navier-Stokes equation~\eqref{ns} with data~$\varphi_n$ for each
integer~$n$, and produce a unique solution~$u_n \in E_{a,b}(T_{u_n})$ for
some time~$T_{u_n} < T^*(\varphi_n)$.
Now let us define, for any~$J \geq 0$
$$
r_n^J (x,t) := u_n (x,t) - \left( \sum_{j = 0}^{J} \Lambda_{j,n}
U_j (x,t) + w_n^J (x,t) \right)\,,
$$
where $ \Lambda_{0,n} {U}_{0}
(x,t) := U_0(x,t)$, and we abuse our earlier notation for
dilations/translations to define
$$
 \Lambda_{j,n} U_j(x,t) := \frac{1}{\lambda _{j,n} }
U_j \left( \frac{x - x_{j,n}}{\lambda _{j,n} }, \frac{t}{\lambda
_{j,n} ^2}\right).
$$
To be consistent, from now on we define $\lambda_{0,n} := 1$ and
$x_{0,n} := 0$, and remark that the lifetime of the re-scaled $U_j$ has
become at least $\lambda_{j,n}^2 T_j$. Therefore, the function~$r_n^J(x,\cdot)$ is defined a priori for~$t \in
[0,t_n]$, where
$$
t_n = \min \left( T_{u_n} ; \min_{j \in I} \lambda _{j,n} ^{2}T_j ; \infty
\right) = \min (\tau_n ; T_{u_n})
$$
with the notation of Theorem~\ref{thm:profa}.  Our main goal consists
in proving that~$r_n^J$ is actually defined on~$[0,\tau_n]$ (at least
if~$J$ is large enough),
which will be a consequence of the
perturbation theory for the Navier-Stokes equation. In the process, we
shall obtain the desired uniform limiting property
$$
 \lim_{J\rightarrow \infty} \left(\limsup_{n \rightarrow \infty}
\|r_n^J\|_{E_{p,p}(\tau_n)}\right) = 0\,.
$$
Let us write the equation satisfied by~$r_n^J$.  It turns out to be easier
to write that equation after a re-scaling in space-time. For
convenience, let us re-order the functions~$\Lambda_{j,n} U_j $, for~$0
\leq j \leq j_0$, in such a way that, for some $n_0 = n_0(j_0)$ sufficiently large, we have
\begin{equation}
\label{eq:increasTglobal} \forall n\geq n_0,\qquad j\leq j' \leq j_0 \quad \Longrightarrow \quad \lambda
_{j,n} ^2 T_j^* \: \: \leq \: \: \lambda_{j',n} ^2 T_{j'}^*
\end{equation} (some of these terms may equal infinity), where~$ T_j^*$ is the maximal life span
of~$\phi_j$ (such a reordering is possible on a fixed and finite number of profiles
due to the orthogonality of scales). In particular, with this ordering we have $\tau_n = \l_{0,n}^2 T_0$, and we note that~$\lambda _{j,n} ^2 T_j^* $ is the life span
of~$\Lambda_{j,n} U_j$.

The inverse of our
dilation/translation operator $\Lambda_{j,n}$ is
\begin{equation}
\label{eq:scalop} \Lambda^{-1}_{j,n} f(s,y) :=
\lambda _{j,n} f ( \lambda _{j,n} y + x_{j,n} , \lambda _{j,n} ^2 s ).
\end{equation}
We then define, for every integer~$J$,
$$
\forall j \leq J, \quad U_n^{j,0} := \Lambda^{-1}_{0,n} \Lambda_{j,n} U_j, \quad
R_n^{J,0} := \Lambda^{-1}_{0,n} r_n^J,
$$
$$
W_n^{J,0} := \Lambda^{-1}_{0,n} w_n^J \quad \mbox{and} \quad U_n^{0} := \Lambda^{-1}_{0,n} u_n.
$$
Clearly we have
$$
 R_n^{J,0} (s,y) := U_n^{0} (s,y) - \left( \sum_{j = 0}^{J}
U_n^{j,0}(s,y) + W_n^{J,0}(s,y) \right) ,
$$
and~$ R_n^{J,0}$ (which  for the time being is defined for times~$s$ in~$ [0,t_n^0)$ where~$t_n^0 := \min \{T_0, \l_{0,n}^{-2}T_{u_n} \}$) solves the following system:
\begin{equation}\label{remaindereq} \displaystyle
\left\{\begin{array}{rcl} \displaystyle \partial_s R_n^{J,0} +
\mathbb{P} ( R_n^{J,0} \cdot \nabla R_n^{J,0} ) -  \Delta R_n^{J,0}
+ Q( R_n^{J,0},F_n^{J,0}) & = & G_n^{J,0} \\ \left
. R^{J,0}_{n}\right|_{s=0} & = & 0 ,
\end{array} \right.
\end{equation}
where we recall that~$\mathbb{P}:=Id-\nabla \Delta^{-1}(\nabla\cdot)$ is the
projection onto divergence free vector fields,
$$
Q(a,b):=\mathbb{P}(a\cdot \nabla b+b\cdot \nabla a)
$$ for two vector
fields $a,b$, and, finally, where
\begin{equation}
  \label{eq:defF}
F_n^{J,0} := \sum_{j \leq J } U_n^{j,0} + W_n^{J,0} ,
\end{equation}
and
\begin{multline}
  \label{eq:defG}
  G_n^{J,0} := - \frac{1}{2}\! \! \! \!  \sumetage{j \neq j' }{(j, j')
\in \{0,..,J\}^2} Q( U_n^{j,0} , U_n^{j',0})\\
{}- \sum_{j \leq J} Q(
U_n^{j,0} , W_n^{J,0} )- \frac{1}{2}Q ( W_n^{J,0}, W_n^{J,0}).
\end{multline}
Notice here that in re-ordering the
profiles, we may have lost the fact that~$\lambda_{0,n}=1$; however we have a
(simpler notation for a) uniform lower bound on the lifetime of all possibly blowing-up profiles: for any real number~$T_0$ smaller than~$T_0^*$, \eqref{eq:increasTglobal} gives
$$
\forall j \in \mathbb{N} , \quad \{U_n^{j , 0}\}_{n\geq n_0} \quad \mbox{is bounded
in} \: E_{a,b}({T_0})\ .
$$
However, we have no uniform control over the sum $F^{J,0}_n$
 which enters the drift term in the
perturbed equation \eqref{remaindereq}. In order to use perturbative
bounds on this system, as stated in Proposition~\ref{pdr1} in the appendix, we need such a
good control on the drift term (which will come from orthogonality arguments), and smallness of the forcing term
$G^{J,0}_n$.

We start with the drift term. Notice that we do not claim uniform
 boundedness in $E_{p,p}$ but rather in a weaker space, which will be
 enough for our purposes (the reader might notice that we could
 replace $p$ in the following statement by $b$ ($<p$) but this will not be necessary and would make notations even more cumbersome).
\begin{lemma} \label{lemmadrift}
The sequence~$(F_n^{J,0})$ is bounded in~$
\mathcal{L}^p([0,T_0];\dot B^{s_p+\frac2p}_{p,p})$, uniformly in~$J$.
\end{lemma}
Let us prove this lemma. Notice that
$$
F_n^{J,0} = \sum_{j \leq j_0} U_n^{j,0} + \sum_{j_0+1}^J U_n^{j,0} +
W_n^{J,0}
$$
and by the small data theory, $U_n^{j,0} $ (for~$j > j_0$)
and~$W_n^{J,0}$ are bounded in~$ E_{p,p}({\mathbb{R}^+})$ by their
respective initial data in~${\dot B^{s_{p}}_{p,p}}$; in particular,
for all~$1\leq r\leq +\infty$,
\begin{equation}\label{boundlastprofiles} \forall j > j_0, \quad
\|U_n^{j,0} \|_{{\mathcal{L}}^{r}(\mathbb{R}^+; \dot B^{s_p + \frac
2r}_{p,p})} = \| {U}_{j }\|_{\mathcal
{L}^{r}(\mathbb{R}^+; \dot B^{s_p + \frac 2r}_{p,p})} \lesssim \|
\phi_j\|_{\dot B^{s_{p}}_{p,p}}.
\end{equation}
Therefore, \eqref{orth3a} yields
$$
\sum_{j=j_0+1}^J \|U^{j,0}_n\|^p_{\mathcal{L}^p_t \dot B^{s_p+\frac2p}_{p,p}}\leq \sum_{j=j_0+1}^J \|U^{j,0}_n\|^p_{E_{p,p}(+\infty)} <+\infty,
$$
 where the bound is uniform in $J$.

We then need to use the orthogonality of cores/scales in the form of the following lemma:
\begin{lemma}
\label{lastmin}  Fix any $r\in [p,\infty)$.  There exists some ${\e(J,n):\N \times \N
\to \R^+}$ satisfying, for each fixed $J$,
$$
\lim_{n\to\infty}\e(J,n) = 0
$$
 and such that, for any
$J'< J$,
$$
\left \|\sum_{j=J'}^{J} U^{j,0}_n
 \right \|^p_{\mathcal{L}^r_t \dot B^{s_p+\frac2r}_{p,p}} \leq   \sum_{j=J'}^{J}
  \|U^{j,0}_n\|_{\mathcal{L}^r_t \dot B^{s_p+\frac2r}_{p,p}}^p + \epsilon(J,n).
 $$
\end{lemma}
Postponing the proof of this for a moment, let us use it to complete the proof of Lemma~\ref{lemmadrift}.
 Lemma~\ref{lastmin}, along with the small data theory for Navier-Stokes, implies that
 \begin{equation}
   \label{eq:sumj0J}
2^{-p} \Big\|\sum_{j=j_0+1}^{J} U_n^{j,0} + W_n^{J,0}
\Big\|^p_{\mathcal{L}^p_t \dot B^{s_p+\frac2p}_{p,p}}\lesssim \sum_{j=j_0+1}^\infty
\|\phi_j\|^p_{B^{s_p}_{p,p}}+\e(J,n)+ \|\psi_n^{J}\|^p_{B^{s_p}_{p,p}}.
 \end{equation}
  The first~$j_0+1$ terms
are then dealt with because the scaling we chose is such that the norm of~$ U_n^{j,0} $
in~$ E_{p,p}({[0,T_0]})$ is bounded uniformly in~$n$, by a constant
depending on~$T_0$, and that concludes the proof of the bound
on~$(F_n^{J,0})$:  Lemma~\ref{lemmadrift} is proved. \qed
\\

Now let us prove Lemma \ref{lastmin}.  Recall that for any $\sigma \in \mathbb{R}$, we have the following equivalence
of norms, where we replace the $2^j$ scale by a continuous parameter~$\tau\sim 2^{-2j}$ (which is easier to handle with rescaling) and the
frequency localization $\Delta_j$ by the derivative of heat kernel~$K(\tau):=\tau \partial_\tau
e^{\tau\Delta}$:
\begin{equation}\label{besequiv} \|f\|_{\dot B^{\sigma}_{p,q}} \sim
\left\| \|\tau^{-\sigma/2} K(\tau) f\|_{L^p}\right\|_{L^q(\mathbb{R}^+,
\frac{d\tau}{\tau})}.
\end{equation}
 There is no difficulty in adding the time norm, and hence, setting
 $$
 \gamma:=-1-p s_p/2-p/r
 $$
 with $1\leq r\leq +\infty$, we have
\begin{equation}\label{defbespp}
 \|f\|_{\mathcal{L}^r_t(\dot B^{s_p+\frac 2 r}_{p,p})}^p \sim
 \int_0^\infty \tau^{\gamma}\|K(\tau)
f\|_{L^r_t L^p_x}^p\ d\tau\ .
\end{equation}
We proceed now with the lemma.  Writing
$$
\left\|K(\tau)\sum_\ell f_ \ell \right\|^p_{L^r_tL^p_x}= \left\|\ \bigl|\sum_ \ell K(\tau) f_ \ell \bigr|^p\
\right\|_{L^\frac r p_t L^1_x}\,,
$$
and applying the elementary inequality
\begin{equation}
\label{ineq}
 \left| \Bigl| \sum_{\ell =1}^L A_ \ell \Bigr|^p
- \sum_{\ell =1}^L |A_ \ell |^p \right| \leq \ \ C(L) \sum_{\ell\neq \ell'} |A_ \ell |
|A_{\ell'}|^{p-1}\ ,
\end{equation}
to the sum inside the norm on the right which we take to be~$\displaystyle \sum_{j=j_0+1}^J K(\tau)U^{j,0}_n$ and applying the triangle inequality in
$L^\frac r p_tL^1_x$ (this is where we need the restriction $r\geq p$) we get $C(J)$ cross-terms like
\begin{multline*}
  \int_0^\infty \tau^{\gamma} \left\| K(\tau) \left[
\frac{1}{\lambda _{1,n}} U_1\left(\frac{\cdot-x_{1,n}}{\lambda
_{1,n}}, t\right)\right](x)\right.\\
{} \times \left.\left(K(\tau) \left[
\frac{1}{\lambda _{2,n}} U_2\left(\frac{\cdot -x_{2,n}}{\lambda
_{2,n}}, t\right)\right](x)\right)^{p-1}\right\|_{L^{\frac r p}_t L^1_x} d\tau =\\
\int_0^\infty \tau^{\gamma}  \left\| \frac{1}{\lambda
_{1,n}}\left[ K({\lambda
  _{1,n}^{-2}\tau})[U_1(t)]\right]\left(\frac{x-x_{1,n}}{\lambda
_{1,n}}\right)\right.\\
{}\times\left. \left( \frac{1}{\lambda
_{2,n}}\left[ K(\lambda _{2,n}^{-2}\tau)
[U_2(t)]\right]\left(\frac{x-x_{2,n}}{\lambda
_{2,n}}\right)\right)^{p-1}\right\|_{L^\frac r p_t L^1_x} d\tau .
\end{multline*}
Using the change of variables
$\displaystyle y=\frac{x-x_{1,n}}{\lambda _{1,n}}$, $s= \lambda _{1,n}^{-2}\tau$, we
see that this term equals
$$\left(\frac{\lambda _{1,n}}{\lambda _{2,n}}\right)^\frac{2+d}{p'}
\int_0^\infty \left\|V_1(y,s) \left(V_2\left(\frac{\lambda
_{1,n}}{\lambda _{2,n}}y + \frac{x_{1,n}-x_{2,n}}{\lambda _{2,n}},
\left(\frac{\lambda _{1,n}}{\lambda _{2,n}}\right)^2s
\right)\right)^{p-1} \right\|_{L^\frac r p_t L^1_y}\ ds ,
$$
where $\frac{1}{p} + \frac{1}{p'} = 1$ and $V_i(y,s) =
s^{\frac{\gamma}{p}}\left[K({s})[U_i(t)]\right](y)$
for $i\in\{1,2\}$.  Since~$U_i(t) \in E_{p,p}$, we can approximate~$V_1$ and~$V_2$ in~$L^p_s L^r_t L^p_x$ by smooth functions
 of~$(s,t,y)$
with compact support in~$(0,+\infty)^2\times \mathbb{R}^d$.  It is
therefore clear (by dislocation of the inner~$L^1_y$) that the above term tends to zero as~$n\to\infty$ if~$\frac{\lambda _{1,n}}{\lambda _{2,n}} \to 0$ as~$n\to\infty$, or if~$\lambda _{1,n}\equiv \lambda _{2,n}$ and~$ \left|\frac{x_{1,n}-x_{2,n}}{\lambda _{2,n}} \right|\to \infty$ as~$n\to\infty$.
Similarly, if we appropriately chose the new variables in terms of~$\lambda _{2,n}$ and~$x_{2,n}$ instead, we can easily show that the
term also tends to zero if~$\frac{\lambda _{2,n}}{\lambda _{1,n}} \to
0$, and the result follows in view of~\eqref{orthseq}. Lemma~\ref{lastmin} is proved. \qed
\\

We now turn to the source term and prove another lemma.
\begin{lemma}\label{lemmaremainder}
Let $G_n^{J,0}$ be the source term defined by \eqref{eq:defG} and
$$
F_{p,p}(T_0):= {\mathcal{L}}^{\frac{2p}{p+1}}([0,T_0];\dot
B^{s_p-1+\frac1p}_{p,p})+{\mathcal{L}}^{p'}([0,T_0];\dot
B^{s_p-\frac2p}_{p,p}),
$$
where~$1/p+1/p' = 1$.
Then
$$
\lim_{J \rightarrow \infty} \limsup_{n \rightarrow \infty} \|
G_n^{J,0}\|_{F_{p,p}(T_0)} = 0.
$$
\end{lemma}
First by standard product laws in Besov spaces
(joint with a H\"older estimate in time) we gather
that
\begin{eqnarray*}
\|Q ( W_n^{J,0}, W_n^{J,0})\|_{\mathcal{L}^{p'}([0,T_0];\dot B^{s_p-\frac2p}_{p,p})} &\leq & \| W_n^{J,0}
\otimes W_n^{J,0}\|_{{\mathcal{L}}^{p'}([0,T_0];\dot B^{s_p+1-\frac2p }_{p,p})} \\ & \lesssim & \|W_n^{J,0}\|_{\mathcal{L}^{2p'}([0,T_0];\dot B^{s_p +\frac1{p'}}_{p,p})}^2.
\end{eqnarray*} Note that the condition~$0<s_p+1/p'=(d-1)/p<d/p$ allows to
justify the product law.

Then by definition of~$ W_n^{J,0} $ and due to the scaling invariance
of the space~${\mathcal{L}}^{2p'}([0,T_0];\dot B^{s_p + \frac 1 {p'}}_{p,p})$ we
obviously have
\begin{equation}\label{limQWW} \lim_{J \rightarrow \infty} \limsup_{n
\rightarrow \infty} \|Q ( W_n^{J,0}, W_n^{J,0})\|_{\mathcal
{L}^{p'}([0,T_0];\dot B^{s_p -  2/ {p}}_{p,p})} = 0.
\end{equation}
Define $fg=\mathcal{T}_f g+\mathcal{T}_g f+\Pi(f,g)$ the paraproduct decomposition of the
product $fg$, and $\zeta(f,g) : =fg-\mathcal{T}_f g$.
We notice that  (abusing
notations as the $Q$ entries are vectors)
$$
\|Q(f,g)\|_{\mathcal
{L}^{p'}([0,T_0];\dot B^{-2 +s_p + \frac 2{p'}}_{p,p})} \lesssim  \|fg\|_{\mathcal
{L}^{p'}([0,T_0];\dot B^{-1 +s_p + \frac 2{p'}}_{p,p})},
$$
and we have
\begin{multline*}
\| \zeta( \sum_{j<J} U_n^{j,0} , W_n^{J,0}) \|_{\mathcal
{L}^{p'}([0,T_0]; \dot B^{-1 +s_p + \frac 2{p'}}_{p,p})} \lesssim \|W_n^{J,0}\|_{{\mathcal{L}}^{r}([0,T_0];\dot B^{s_p + \frac
2 r}_{p,p})} \\
{}\times \Bigl\|\sum_{j<J}  U_n^{j,0} \Bigr\|_{{\mathcal{L}}^{p}([0,T_0];\dot B^{s_p+\frac2p}_{p,p})}
\end{multline*}
where $1/p'=1/r+1/p$.  The product law is allowed
because the sum of the two respective
regularities is again $2s_p+2/r+2/p=2(d-1)/p>0$. Notice that the
regularity of~$W^{J,0}_n$   is~$\sigma=d/p + (p-4)/p > d/p$ for large~$p$, which
explains why we subtracted the paraproduct term carrying the
corresponding low frequencies of $f$. This specific term is handled
differently, as
\begin{multline*}
\| \mathcal{T}_{ (\sum_{j<J} U_n^{j,0})}  W_n^{J,0} \|_{\mathcal
{L}^{\frac{2p}{p+1}}([0,T_0];\dot B^{s_p -1 + \frac{p+1}{p}}_{p,p})} \lesssim \|W_n^{J,0}\|_{{\mathcal{L}}^{\frac{2p}{p-1}}([0,T_0];\dot B^{s_p + \frac
{p-1} p}_{p,p})}\\
{}\times \Bigl\|\sum_{j<J}  U_n^{j,0} \Bigr\|_{{\mathcal{L}}^{p}([0,T_0];\dot B^{s_p+\frac2p}_{p,p})}
\end{multline*}
where we do not care whether the sum of regularities $-1+(2d+1)/p$ is
negative, thanks to the frequency distribution in the paraproduct.

The profiles $U_n^{j,0}$ in~${\mathcal{L}}^{p}([0,T_0];\dot B^{s_p+\frac2p
}_{p,p})$ are obviously under control for~$j > j_0$, as seen
in the proof of Lemma~\ref{lemmadrift} (see~ \eqref{boundlastprofiles}), and so is their sum due to Lemma~\ref{lastmin}.    It follows that
\begin{equation}\label{limQUW}
\lim_{J \rightarrow \infty} \limsup_{n\rightarrow \infty} \:
\bigl\|Q (  \sum_{j<J} U_n^{j,0},
W_n^{J,0})\bigr\|_{F_{p,p}(T_0)} =0.
\end{equation}
Now we are left with the terms~$ Q( U_n^{j,0} ,
U_n^{j',0})$ for~$j \neq j'$. Again, we resort to the
orthogonality of the profiles to conclude: here we use the fact that,
at fixed $j<J$, $U^{j,0}_n$ belongs to~$E_{p,p}(T_0)$, but with no uniform bound
with respect to $J$; by scaling we have (with $r$ to be specified later)
$$
\bigl \| Q( U_n^{j,0} , U_n^{j',0}) \bigr\|_{{\mathcal{L}}^r_{T_0} (\dot B^{-2
+ s_p + \frac 2r}_{p,p})} =\bigl \| Q( \Lambda_{j,n}U_j , \Lambda_{j',n}U_{j'}) \bigr\|_{{\mathcal{L}}^r_{(\lambda_n^0)^2T_0}(\dot B^{-2 +s_p +
\frac 2r}_{p,p})}
$$
and we recall that
$$
 \Lambda_{j,n} U_j(x,t): = \frac{1}{\lambda _{j,n} } U_j \left( \frac{x
- x_{j,n}}{\lambda _{j,n} }, \frac{t}{\lambda _{j,n} ^2}\right)
$$
is defined at least on~$[0,\lambda_{0,n}^2T_0]$ by construction.  \\\\
Suppose to simplify that~$T_j \leq T_{j'}$ (if not exchange~$j$
and~$j'$).  By density (for $r<\infty$) for any~$\varepsilon >0$ one can find smooth
compactly supported functions~$F^j$ and~$F^{j'}$ such that
$$
\|F^j - U_j \|_{{\mathcal{L}}^{2r}([0,T_j];\dot B^{s_p + \frac
1r}_{p,p})} + \|F^{j'} - U^{j'} \|_{{\mathcal{L}}^{2r}([0,T_j];\dot
B^{s_p + \frac 1r}_{p,p})} \leq \varepsilon.
$$
Then we decompose
$$
Q( U_n^{j,0}  , U_n^{j',0} ) = Q( \rho_n^{j} , \rho_n^{j'})
+ Q( f_n^{j} , \rho_n^{j'}) + Q( \rho_n^{j} , f_n^{j'})+ Q( f_n^{j}
,f_n^{j'})
$$
where
$$
 f_n^{j} := \frac{1}{\lambda _{j,n} } F^j \left( \frac{x -
x_{j,n}}{\lambda _{j,n} }, \frac{t}{\lambda _{j,n} ^2}\right)
$$
and~$U_n^{j,0}  = f_n^{j} + \rho_n^j$.  The same argument as
above enables us to write that for all~$j, j'$,
\begin{equation}\label{limQrhorho} \sup_{n} \|Q(
\rho_n^{j} , \rho_n^{j'}) \|_{{\mathcal{L}}^r([0,T_0];\dot B^{-2+s_p +
\frac 2r}_{p,p})} \leq C\varepsilon^2
\end{equation} and
\begin{equation}\label{limQrhofbis} \sup_{n } \|Q(
f_n^{j} , \rho_n^{j'}) +Q( \rho_n^{j} , f_n^{j'}) \|_{\mathcal
{L}^r([0,T_0];\dot B^{-2 +s_p + \frac 2r}_{p,p})} \leq C \varepsilon,
\end{equation} where~$C$ only depends on~$T_0$ (neither on~$j$ nor
on~$j'$). So finally we are left with~$ Q( f_n^j ,f_n^{j'})$, or more
accurately with the product $f_n^j f_n^{j'}$. But that product is
dealt with exactly as the cross terms in the proof of Lemma
\ref{lastmin}, as smoothness lets us estimate the product in the space~$\mathcal
L^r([0,T_0]; \dot B^{s_r+\frac2r-1}_{r,r})$ or as in Proposition~2.1
of~\cite{Gallagher}, up to a suitable modification of the time
exponents to reach $1<r$.

Chosing $r=p'$, by Sobolev's embedding this implies in particular that
\begin{equation}
  \label{limQff}
       \lim_{n \rightarrow \infty} \|Q( f_n^j , f_n^{j'})
\|_{\mathcal{ L}^{p'}([0,T_0];\dot B^{  s_p - \frac 2 p}_{p,p})} = 0.
\end{equation}
Putting together estimates~\eqref{limQWW}
to~\eqref{limQff} ends
the proof of   Lemma~\ref{lemmaremainder}.\qed

Using Lemma~\ref{lemmadrift}, we have that the drift term $F^{J,0}_n$ is
uniformly bounded
in $\mathcal{L}^p([0,T_0];\dot B^{s_p+\frac2p}_{p,p})$ by a constant $C_0$
(depending on the profiles, on~$j_0$ and on~$T_0$), provided
that at fixed $J$, $n$ is chosen larger than some $N(J)$. Then if one
picks $J$ large enough, the forcing term is small enough (with respect
to the smallness condition in the perturbation result recalled in Proposition~\ref{pdr1}), and we obtain the desired
bound
on $R^{J,0}_n$ in $E_{p,p}(T_0)$ for $J>J_0$ and $n>N(J_0)$ thanks to
Proposition~\ref{pdr1}. This bound immediately implies that $u_n$ is
bounded in $E_{p,p}(T_0)$, and by persistence of regularity (see
e.g. \cite{gip3}), we get $u_n\in E_{a,b}(T_0)$.

\medskip

All other statements
in Theorem  \ref{thm:profa} follow easily, as in \cite{Gallagher},
which ends the proof. \qed
\section{Serrin's endpoint regularity criterion}
\label{serrin}
\subsection{Preliminaries and statement of the main result}
   Let $X=X(\rd)$ be a critical space with respect to the
Navier-Stokes scaling, that is
$$
\|\lambda f(\lambda \cdot)\|_X =
\|f\|_X
$$ for any $f\in X$, and such that local in time strong solutions
can be constructed.  Let $T^*(u_0)$ be the maximal time of existence
of the solution~$u=NS(u_0)$   in~$X$ for any $u_0\in X$.  The type of statement we
would like to address is the following:
\begin{statement}\label{statementa} For any $u_0 \in X$,
$$\sup_{t\in [0,T^*(u_0))}\|NS(u_0)(t)\|_X < \infty \qquad  \Longrightarrow \qquad T^*(u_0) = +\infty\ .$$
\end{statement}
Such a statement was proved in \cite{ess} for
$X={L^3(\mathbb{R}^3)}$ (in the context of Leray-Hopf weak solutions).
An alternative proof of that theorem was given in \cite{kk} for $X =
{\dot H^\frac{1}{2}}(\mathbb{R}^3)$ using the profile decompositions
in~\cite{Gallagher}.  (Of course that result is implied as well by
\cite{ess}.)  Our goal now is to give a proof similar to the one in \cite{kk}
using the profile decompositions in this article for the case $X=
{L^3(\mathbb{R}^3)}$. This will give a different proof of the
following, which was proved in \cite{ess} and also extended to $d>3$ in~\cite{dongdu2}:
\begin{thm}[Endpoint regularity criterion]\label{regularity}
For any $u_0 \in {L^3(\mathbb{R}^3)}$,
$$\sup_{t\in [0,T^*(u_0))}\|NS(u_0)(t)\|_{L^3(\mathbb{R}^3)} < \infty \qquad  \Longrightarrow \qquad T^*(u_0) = +\infty\ .$$
\end{thm}
   Note that due to the time-continuity in
${L^3(\mathbb{R}^3)}$ of strong solutions, the left-hand side is
equivalent to $NS(u_0) \in
L^\infty((0,T^*(u_0));{L^3(\mathbb{R}^3)})$, or in the notation of \cite{ess},  $NS(u_0) \in L_{3,\infty}(0,T^*(u_0))$.
\subsection{Proof of Theorem \ref{regularity} }
   Theorem \ref{regularity}   follows from Theorems~\ref{thm:aa} to~\ref{thm:ac} below, which will be proved in the
subsequent sections.  In the following, we define $A_c>0$ by
$$A_c := \sup \left\{A>0\ ;\ \sup_{t\in
[0,T^*(u_0))}\|NS(u_0)(t)\|_{L^3(\mathbb{R}^3)} \leq A \Rightarrow
T^*(u_0) = \infty \right\}.
$$
Note that $A_c$ is well-defined by small-data results.  Moreover,
  if $A_c $ is finite, then we have
$$
A_c = \inf \left\{\sup_{t\in [0,T^*(u_0))}\|NS(u_0)(t)\|_{L^3(\mathbb{R}^3)} \ ; \  T^*(u_0) < \infty \right\}\ .
$$
\begin{thm}[Existence of a critical element]\label{thm:aa}
 Suppose $A_c < +\infty$, and let $u_{0,n}$
be any bounded sequence in~$L^3$ such that $T^*(u_{0,n}) < +\infty$ and such that~$ A_c \leq
A_n$, where
$$
A_n:= \sup_{t\in [0,T^*(u_{0,n}))} \|NS(u_{0,n})(t)\|_{L^3} <\infty\ ,$$
and satisfying $A_n \to A_c$ as $n\to\infty$.  Let $U_j = NS(\phi_j)$ be the
Navier-Stokes profiles associated to~$\{u_{0,n}\}$.  Then there exists an integer~$j_0 \in \mathbb{N} $
such that
$$T^*(\phi_{j_0 }) < +\infty \quad \mbox{and} \quad  \sup_{t\in
[0,T^*(\phi_{j_0 }))} \|U_{j_0}(t)\|_{L^3} = A_c.
$$
\end{thm}
   We shall call any solution with the properties of $U_{j_0}$
above a ``critical element" --- that is, those solutions whose supremum in $L^3$
is $A_c$ and who blow up in finite time.  For the next two statements we fix~$u_0$ in~$L^3$.
The first gives a
kind of compactness property for critical elements:
\begin{thm}[Compactness of   critical elements]\label{thm:ab}
Suppose that~$A_c  $ is finite and that~$T^*(u_{0}) < +\infty$ and assume
$$
\displaystyle \sup_{t\in
[0,T^*(u_{0}))} \|u(t)\|_{L^3} = A_c
$$
where $u=NS(u_0)$. Then $u(t)
\to 0$ in $\mathscr{S}'$ as $t \nearrow T^*(u_0)$.
\end{thm}
Finally, we state the following, which corresponds to the ``rigidity theorem" in \cite{kk}:
\begin{thm}[Global existence of critical elements]\label{thm:ac}
Assume that the solution~$u:=NS(u_0)$ satisfies
$$
\displaystyle \sup_{t\in [0,T^*(u_{0}))} \|u(t)\|_{L^3} <+\infty$$
and moreover $u(t) \to 0$ in $\mathscr{S}'$ as $t
\nearrow T^*(u_0)$. Then $T^*(u_0) = +\infty$.
\end{thm}
   Theorems \ref{thm:aa} to~\ref{thm:ac} immediately imply that
$A_c = +\infty$, and Theorem~\ref{regularity} is proved.  Indeed, if
$A_c < +\infty$ then Theorem \ref{thm:aa} guarantees us some~$j_0$
such that $u_0 = \phi_{j_0 }$ satisfies the conditions of Theorem
\ref{thm:ab}, but then Theorem \ref{thm:ac} contradicts the fact that
$T^*(u_0) < +\infty$.  \hfill $\Box$

In order to prove Theorems
\ref{thm:aa} to~\ref{thm:ac}, we shall need the following result, which
was proved in \cite{gip3}:
\begin{thm}[\cite{gip3}]\label{decay} Let $u_0$ be
some divergence-free vector field in~$ X$ such that $T^*(u_0) = +\infty$, where
$X$ is either ${L^3(\mathbb{R}^3)}$ or ${\dot
B^{s_p}_{p,p}}(\mathbb{R}^3)$.  Then $$\displaystyle{\lim_{t\to
+\infty}\|NS(u_0)(t)\|_X = 0}\ .$$
\end{thm}
   We shall also need the following lemma in the spirit of Lemma \ref{lastmin} whose proof we outline
below (see the proof of (3.18) in \cite{kk} for more details):
\begin{lemma}\label{clm:orth} Suppose $\displaystyle \sup_n
\|u_{0,n}\|_{L^3(\mathbb{R}^3)} < \infty$ and let $\{U_j\}$ be the
associated Navier-Stokes profiles given by Theorem \ref{thm:profa}.
Let $\{t_n\} \subset \mathbb{R}^+$ be any sequence such that $t_n \leq
\tau_n$ for all $n$. There exists a subsequence in~$n$ such that the
following is true: for any $\epsilon >0$ and $J\in \mathbb{N} $, there
exists~$N_0 = N_0(J,\epsilon ) \in \mathbb{N} $ such that, for any $J'
< J$,
$$\left| \bigl\|\sum_{j=J'}^{J} \Lambda_{j,n}U_j(t_n)\bigr\|_\lt^3 - \sum_{j=J'}^{J} \|\Lambda_{j,n}U_j(t_n)\|_\lt^3\right| < \epsilon  $$
for all $n\geq N_0$.
\end{lemma}
To prove this lemma, using again the elementary inequality \eqref{ineq}, we must estimate a finite number (depending on $J$) of
terms of the form
$$
\int |\Lambda_{j_1,n}U_{j_1}(t_n) |\ |\Lambda_{j_2,n}U_{j_2}(t_n)|^{2}
$$
 where $j_1,j_2 \in \{J',\dots,J\}$ and $j_1 \neq j_2$.
 Such a term
tends to zero uniformly (for each fixed $J$) as $n\to\infty$ after
passing to the following subsequence: if $t_n \lambda _{j_i,n}^{-2}
\to \infty$ for $i=1$ or $i=2$, then necessarily the $i$th profile is
global and using H\" older once more this term tends to zero by
Theorem \ref{decay}.  If neither tends to infinity, we may pass to a
subsequence so that the re-scaled times in both terms approach a
constant time. Since we stay strictly away from the blow-up times of
any profile, we may use the time continuity of solutions to replace
the re-scaled time value by this fixed time in each profile, after
which the orthogonality of the scales/cores implies that the term
tends to zero due to \eqref{Ldorthogx}.  One may take a diagonal subsequence so that this is
true in all instances when the times are bounded which proves Lemma~\ref{clm:orth}.  \qed
\subsection{Proof of Theorem \ref{thm:aa}} Consider the bounded
sequence described in Theorem \ref{thm:aa} and its profile
decomposition after passing to a subsequence.  For notational
convenience, set $\lambda _{0,n} \equiv 1$ and $x_{0,n} \equiv 0$, so
that the ``0'th" profile is the weak limit (without any
transformations).

Note first that since $T^*(u_{0,n}) < +\infty$
for all $n$ in the sequence which we consider, there must be at least
one $j\geq 0$ such that $T^*(\phi_j) < \infty$.  If not, one could
take $\tau_n \equiv \infty$ in Theorem \ref{thm:profa} and hence $u_n$ is global for large $n$, contrary to assumption.

Property \eqref{orth3} and small
data results in ${L^3(\mathbb{R}^3)}$ now imply that there exists some
$J_0 \in \mathbb{N}  $ such that after re-ordering the
profiles one has $$T_j^* < \infty \iff 0\leq j \leq J_0. $$
The
orthogonality properties of the scales $\lambda _{j,n} $ then allow us
to re-arrange the first $J_0$ profiles in such a way that  for all sufficiently large~$n$, one has~$\lambda
_{0,n}^2 T^*_0 \leq \lambda _{1,n}^2 T^*_1 \leq \cdots \leq \lambda
_{J_0,n}^2 T^*_{J_0}$.  Fix now such an
ordering of the $\phi_j$ (so it is now possible that $\lambda _{0,n}
\neq 1$, $x_{0,n} \neq 0$).

 Note now that for any $s \in
[0,T_0^*)$,   we have
\begin{equation}
\label{zza}
 \Lambda^{-1}_{0,n} u_n ( s) =
U_0(s)
{} + \sum_{j=1}^J \Lambda^{-1}_{0,n} \Lambda_{j,n} U_j(s)
  +\Lambda^{-1}_{0,n} w_n^J (s) +\Lambda^{-1}_{0,n} r_n^J (s)\ .
\end{equation}
We now claim that the left-hand side converges (after
possibly passing to a subsequence) weakly in ${L^3(\mathbb{R}^3)}$ to
$U_0(s)$, in which case by properties of weak limits and the invariance of the spatial norm we have
$$\|U_0(s)\|_{L^3(\mathbb{R}^3)} \leq \liminf_{n\to\infty} \|u_n(\lambda _{0,n}^2s)\|_{L^3(\mathbb{R}^3)} \leq A_c$$
hence
$$\sup_{s\in [0,T^*_0)} \|U_0(s)\|_{L^3(\mathbb{R}^3)} \leq A_c\ .$$
On the other hand, since $T^*_0 < \infty$, by definition of $A_c$ we
must also have $\sup_{s\in [0,T^*_0)} \|U_0(s)\|_{L^3(\mathbb{R}^3)}
\geq A_c$ and hence $U_0$ is a critical element.

We shall now prove this weak
convergence using the smallness of the
remainders, the time-continuity of the evolution of the profiles and
the decay of the global ones as well as the orthogonality of the
scales/cores.  To simplify notation, in the following we shall denote $\|f\|_p:=\|f\|_{L^p_x}$.

Fix $\epsilon >0$.  We need to show that there exists a subsequence in $n$ such that for any $\varphi \in
\mathcal{C}^\infty_0(\mathbb{R}^3)$,
\begin{equation}\label{wkzero} |<\lambda _{0,n} u_n (\lambda _{0,n}
\cdot + x_{0,n}, \lambda _{0,n}^2 s)-U_0(s), \varphi > | < \epsilon
\end{equation} for $n$ sufficiently large, where
$<\cdot,\cdot>$ denotes the pairing between $L^3$ and $(L^3)' =
L^{\frac{3}{2}}$, i.e. integration over $\mathbb{R}^3$ of the
product.  This will be accomplished by estimating the left-hand side
of~\eqref{wkzero} by
\begin{multline}\label{zzc} \sum_{j=1}^{J_1}
\big|<\Lambda_{j,n}\Lambda^{-1}_{0,n}U_j(s),\varphi> \big| + \left\|\sum_{j=J_1+1}^J
\Lambda_{j,n}\Lambda^{-1}_{0,n}U_j(s)\right\|_3 \|\varphi\|_{\frac{3}{2}}\\ {}+
<\Lambda^{-1}_{0,n}{w}_n^J(s) + \Lambda^{-1}_{0,n}{r}^J_n(s),\varphi>
\end{multline}
where $J_1$ will be some fixed appropriately large integer.  According to~\eqref{orth3}, $J_1$ can be chosen so large that for any $j>J_1$, one has
$${\|U_j(0)\|_3 =\|\phi_j\|_3 \leq \varepsilon_0}$$ and $U_j = NS(\phi_j)$ can be
produced by a fixed-point argument on $(0,\infty)$ (see, e.g.,
\cite{gip3}).  Moreover,
\begin{equation}\label{ldenergy} \sup_{t\geq 0} \|NS(\phi_j)(t)\|_3
\lesssim \|\phi_j\|_3
\end{equation} as a by-product of the small data theory. Hence, for such a $J_1$ we have
\begin{equation}\label{zzdprime} \sum_{j=J_1+1}^J \left\|
\Lambda_{j,n}\Lambda^{-1}_{0,n}U_j(s)\right\|_3^3 \leq C_0 \sum_{j=J_1+1}^J \left\| \phi_j
\right\|_3^3
\end{equation} for some universal $C_0 >0$.
Moreover, by up to a harmless rescaling by $\Lambda^{-1}_{0,n}$ in
its statement (as it does not change the orthogonality of scales/cores), we may use Lemma \ref{clm:orth} and pass to a subsequence in~$n$ so that for each $J$ and any $J' < J$, there exists $n_0(J)$ (in
fact, independent of $J'$) such that for $n\geq n_0(J)$ one has
\begin{equation}\label{zzd}
 \left\|\sum_{j=J'}^J
\Lambda_{j,n}\Lambda^{-1}_{0,n}U_j(s)\right\|_3^3 \leq 2 \sum_{j=J'}^J\left\|
\Lambda_{j,n}\Lambda^{-1}_{0,n}U_j(s)\right\|_3^3\ .
\end{equation}
Now, due to \eqref{orth3}, by taking $J_1$ large enough
(depending on $\varphi$), setting~$J' = J_1+1$ in \eqref{zzd} and
using the previous estimate we can make the middle term of \eqref{zzc}
less than ${\epsilon /3}$ whenever $n \geq n_0(J)$ for $J> J_1$.   Now, by property \eqref{remainders} and the scaling of the norm, there
exists~$n(J) \geq n_0(J)$, increasing in $J$, such that
$$
\|\Lambda^{-1}_{0,n(J)}{w}_{n(J)}^J(s)\|_{\dot
B^{s_{p}}_{p,p}},\|\Lambda^{-1}_{0,n(J)}{r}_{n(J)}^J(s)\|_{\dot
B^{s_{p}}_{p,p}} \to 0
$$
as $J\to\infty$.  In particular, these limits hold weakly.  By heat
estimates,~\eqref{orth4} and the transformational invariance of the
norm, $\Lambda^{-1}_{0,n}{w}_n^J(s)$ is bounded in $L^3$.  We therefore also have
$\Lambda^{-1}_{0,n}{r}_n^J(s)$ bounded in~$L^3$ by~\eqref{zza}, our assumption that
$u_n(t)$ is uniformly bounded in $L^3$ for all times and~\eqref{zzd}
with $J' =0$ in conjunction with (\ref{orth3}) to bound the other
terms.  Therefore the
error terms tend weakly to zero in $L^3$ as well, and hence, setting $n=n(J)$, the third term in
(\ref{zzc}) can be made less than ${\epsilon /3}$ for sufficiently
large $J$.

Finally, since $n(J) \to \infty$ with~$J $, the
orthogonality of the scales/cores shows that each term in the sum on
the left in (\ref{zzc}) tends to zero after a subsequence, by
arguments similar to those in the proof of Lemma~\ref{clm:orth}.
Since there are only a finite number $J_1$ of these, the first term in
(\ref{zzc}) can be made less than ${\epsilon /3}$ for sufficiently
large $J$ which proves (\ref{wkzero}) for some subsequence of $n$'s as
desired.  \hfill $\Box$
\subsection{Proof of Theorem \ref{thm:ab}} Suppose now $u=NS(u_0)$ is
a critical element, and consider the bounded sequence $u_{0,n}:=
u(s_n)$ for some~$s_n \nearrow T^*(u_0)$.  Pass to a subsequence so
that one may write~$u_{0,n}$ and~$u_n:=NS(u_{0,n})$ in terms of
profiles with the notations of Theorems~\ref{thm:dataprofa} and~\ref{thm:profa}.

 As in the proof of Theorem \ref{thm:aa} above,
there is some~$J_0 \geq 0$ such that~$T^*_j < \infty \iff 0\leq j \leq
J_0$ and we may re-arrange the first~$J_0$ profiles in such a way that
$\lambda _{0,n}^2 T^*_0 \leq \lambda _{1,n}^2 T^*_1 \leq \cdots \leq
\lambda _{J_0,n}^2 T^*_{J_0}$ for all sufficiently large $n$.  Fix now
such an ordering of the $\phi_j$, and suppose that $0$ has been moved
now to some $j_0\in \mathbb{N}  $, that is, now $\lambda
_{j_0,n} \equiv 1$ and~$x_{j_0,n} \equiv 0$ and $\phi_{j_0}$ is the
weak limit of $u_{0,n}$. The theorem will therefore be proved if we show that $\phi_{j_0}=0$.

By the definition of $\tau_n$, etc., in
Theorem \ref{thm:profa}, it is clear that with this ordering we must
have
\begin{equation}\label{smallscale} \lambda _{0,n}^2 T^*_0 \leq
T^*(u_{0,n}) = T^*(u_0) - s_n
\end{equation} for large $n$, and hence $\lambda _{0,n} \to 0$ as
$n\to\infty$.  In particular, we see that~$j_0 \neq 0$, that is, $1$
(the scale of the weak limit profile) cannot be smaller than all other
scales.

We shall need the following crucial claim (which actually applies to the more general sequence $\{u_{0,n}\}$ considered in the proof of Theorem~\ref{thm:aa}), whose proof we
postpone momentarily:
\begin{clm}\label{clm:d} Fix any $s\in (0,T^*_0)$.  Setting $t_n:=
\lambda _{0,n}^2s$, after possibly passing to a subsequence in $n$ one
has
$$
\|u_n(t_n)\|_{L^3}^3 \geq \|\Lambda_{0,n}U_0(t_n)\|_{L^3}^3  +
\|u_n(t_n) - \Lambda_{0,n}U_0(t_n)\|_{L^3}^3 + \circ(1)
$$
as $n\to\infty$.
\end{clm}
   Let us proceed to prove Theorem \ref{thm:ab}.  Exactly as in
the proof of Theorem~\ref{thm:aa}, we see again that $U_0$ is a
critical element since we have
$$A_n:= \sup_{t\in [0,T^*(u_{0,n}))} \|NS(u_{0,n})(t)\|_{L^3} = \sup_{t\in [s_n,T^*(u_{0}))} \|u(t)\|_{L^3} \equiv A_c $$
for all $n$, due to the definition of $A_c$ and the fact that
$T^*(u_0) < \infty$.  We shall now show that this implies by Claim
\ref{clm:d} that $\phi_{j_0} = 0$, i.e. that~$u_{0,n}$ tends weakly to zero
which was our goal.  Fix any $\epsilon >0$.  By the time-continuity of
solutions we may take $s\in (0,T^*_0)$ such that
$$
A_c^3 -
\|U_0(s)\|_3^3 < ({\epsilon /2})^3C_0^{-3}
$$ where $C_0>0$ is the
universal constant in the continuous embedding~$L^3 \hookrightarrow
{\dot B^{s_p}_{p,p}}$. Set $t_n:=\lambda _{0,n}^2s$. Then due to Claim
\ref{clm:d}, after passing to a further subsequence in $n$, we have
$$\longformule{
A_c^3 \geq \|u_n(t_n)\|_3^3 \geq \|U_0(s)\|_3^3  }{{}+ C_0^{-3}\|\sum_{j=1}^{J} \Lambda_{j,n}U_j(t_n)+ w_n^J(t_n) + r_n^J(t_n) \|_{{\dot B^{s_p}_{p,p}}}^3 + C_0^{-3}\epsilon (n,s)}
$$
for any $J$ where $\epsilon (n,s) \to 0$ as $n\to \infty$.  According
to (\ref{remainders}), we may fix $J \geq j_0$ so large that
$$\|w_n^J(t_n) + r_n^J(t_n) \|_{\dot B^{s_p}_{p,p}} \leq {\epsilon /2}$$
for sufficiently large $n$.  The previous two inequalities give
\begin{multline*}
  \left(({\epsilon /2})^3 - \epsilon (n,s)\right)^\frac{1}{3} +
  {\epsilon /2} \geq   \|\sum_{j=1}^{J} \Lambda_{j,n}U_j(t_n)\|_{\dot
    B^{s_p}_{p,p}} \\
=: \left( \sum_{j=1}^{J}
\| \Lambda_{j,n}U_j(t_n)\|_{\dot B^{s_p}_{p,p}}^p - \epsilon _J(n)
\right)^\frac{1}{p}\ ,
\end{multline*}
 and we claim (as in Lemma \ref{clm:orth}) that $\epsilon _J(n) \to 0$ as
$n\to \infty$ for fixed $J$, after passing to a subsequence in $n$.
Postponing this fact for a moment, we have now shown that
$$
 \| U_{j_0}({t_n/\lambda _{j_0,n}^2})\|_{\dot B^{s_p}_{p,p}}^p \leq
 \left(\left(({\epsilon /2})^3 - \epsilon (n,s)\right)^\frac{1}{3} +
   {\epsilon /2}\right)^p + \epsilon _J(n)\ .
$$
Recall that $\lambda _{j_0,n} \equiv 1$ because $\varphi_{j_0}$ is the
weak limit of $u_{0,n}$, and note that $t_n = \lambda _{0,n}^2s \to 0$
as $n\to\infty$ due to (\ref{smallscale}).  Therefore letting~$n\to
\infty$ and using the continuity of $U_{j_0}$ in ${\dot
B^{s_p}_{p,p}}$ at $t=0$ we have
$$
\|U_{j_0}(0)\|_{\dot B^{s_p}_{p,p}}
= \|\varphi_{j_0}\|_{\dot B^{s_p}_{p,p}} \leq \epsilon.
$$
Since
$\epsilon $ was arbitrary and $\phi_{j_0} \in L^3$, this implies that
$\phi_{j_0} = 0$, which proves the theorem.

All that remains now
is to show that $\epsilon _J(n) \to 0$, which we now explain.  It is
again a
simple consequence of the orthogonality of the scales/cores, and is
proved in the same way as Lemmas \ref{clm:orth} and~\ref{lastmin}. In fact, up to undoing a harmless
$\Lambda^{-1}_{0,n}$ transform, we follow closely the proof of Lemma~\ref{lastmin}, without the inner $L^r_t$ norm and with  different
times in the profiles $U_i$: using the
elementary inequality \eqref{ineq}, in the expansion of the Besov norm
of the sum $\displaystyle \sum_{j=j_0+1}^J U^{j,0}_n$, we get $O(J)$ cross-terms like
\begin{multline*}
  \int_0^\infty \tau^{\gamma} \int \left|K(\tau) \left[
\frac{1}{\lambda _{1,n}} U_1\left(\frac{\cdot-x_{1,n}}{\lambda
_{1,n}}, s_1\right)\right](x)\right|\\
{} \times \left|K(\tau) \left[
\frac{1}{\lambda _{2,n}} U_2\left(\frac{\cdot -x_{2,n}}{\lambda
_{2,n}}, s_2\right)\right](x)\right|^{p-1}\ dx\ d\tau =\\
\int_0^\infty \tau^{\gamma} \int \left| \frac{1}{\lambda
_{1,n}}\left[ K({\lambda
  _{1,n}^{-2}\tau})[U_1(s_1)]\right]\left(\frac{x-x_{1,n}}{\lambda
_{1,n}}\right)\right|\\
{}\times \left| \frac{1}{\lambda
_{2,n}}\left[ K(\lambda _{2,n}^{-2}\tau)
[U_2(s_2)]\right]\left(\frac{x-x_{2,n}}{\lambda
_{2,n}}\right)\right|^{p-1}\ dx\ d\tau
\end{multline*}
for some $s_1,s_2 >0$ in the life-spans of $U_1$ and $U_2$
respectively (here we have passed to a subsequence and used the time continuity of the profiles and Theorem \ref{decay} as in the proof of Lemma \ref{clm:orth}).  Using the change of variables
$\displaystyle y=\frac{x-x_{1,n}}{\lambda _{1,n}}$, $s= \lambda _{1,n}^{-2}\tau$, we
see that this term equals
$$\left(\frac{\lambda _{1,n}}{\lambda _{2,n}}\right)^\frac{2+d}{p'}
\int_0^\infty \int \left|V_1(y,s)\right| \left|V_2\left(\frac{\lambda
_{1,n}}{\lambda _{2,n}}y + \frac{x_{1,n}-x_{2,n}}{\lambda _{2,n}},
\left(\frac{\lambda _{1,n}}{\lambda _{2,n}}\right)^2s
\right)\right|^{p-1} \ dy\ ds ,
$$
where $\frac{1}{p} + \frac{1}{p'} = 1$ and $V_i(y,s) =
s^{\frac{\gamma}{p}}\left[K({s})[U_i(s_i)]\right](y)$
for $i\in\{1,2\}$.  Since $U_i(s_i) \in {\dot B^{s_p}_{p,p}}$, by
(\ref{besequiv}) we can approximate $V_1$ and $V_2$ in
$L^p(\mathbb{R}^d \times (0,+\infty))$ by smooth functions of $(y,s)$
with compact support in $\mathbb{R}^d \times (0,+\infty)$.  It is
therefore clear that the above term tends to zero as~$n\to\infty$ if
$\frac{\lambda _{1,n}}{\lambda _{2,n}} \to 0$ as $n\to\infty$, or if
$\lambda _{1,n}\equiv \lambda _{2,n}$ and
$\displaystyle \left|\frac{x_{1,n}-x_{2,n}}{\lambda _{2,n}} \right|\to \infty$ as $n\to\infty$.
Similarly, if we appropriately chose the new variables in terms of
$\lambda _{2,n}$ and $x_{2,n}$ instead, we can easily show that the
term also tends to zero if $\frac{\lambda _{2,n}}{\lambda _{1,n}} \to
0$, and the result follows in view of (\ref{orthseq}).\qed

\begin{remark}\label{onlyoneprofile}  A similar argument can be used to show that only one
profile can be a critical element since all others are small at some
time, implying that they exist globally by the small data theory.
Although this fact was used to prove the theorems in \cite{kk}, we
shall not use it here.
\end{remark}

\begin{remark}
One could also prove a similar compactness result as Theorem 3.2 in~\cite{kk}, namely that if~$NS(u_0)$ satisfies
$$
\displaystyle \sup_{t\in
[0,T^*(u_{0}))} \|NS(u_0) (t)\|_{L^3} = A_c
$$
then for any sequence $\{t_n\}$ converging to $T^*(u_{0})$, there exists a sequence $\{s_n\}$ with $t_n \leq s_n \nearrow T^*(u_0)$ such that  sequence~$NS(u_0) (s_n)$ is compact in~$\dot B^{s_p}_{p,p}$ up to norm-invariant transformations in space.
\end{remark}

\begin{remark}
Claim \ref{clm:d} also immediately proves Theorem \ref{thm:aa}, but we feel that the proof given above is more self-contained and perhaps more natural at that point.
\end{remark}

 {\bf Proof of Claim \ref{clm:d}.} \quad  Note first that we may assume without loss of generality that $u_n$ is
scalar-valued by setting $$\|(f^k)_{k=1}^3\|_{L^3} :=
\left\|(\|f^{k}\|_{L^3})_{k=1}^3 \right\|_{\ell^3}$$ and treating each
component separately.

 We first remark that, after passing to
an appropriate subsequence, \def\supetage#1#2{
\sup_{\scriptstyle {#1}\atop\scriptstyle {#2}} }
\begin{equation}\label{consts}
\begin{array}{c} \displaystyle{C_1:= \supetage{J\geq 0}{ n\geq N_0(J)} \|\sum_{j=0}^J \Lambda_{j,n}U_j(t_n)\|_3 \quad < \quad
\infty \quad \quad \textrm{and}} \\\\ \displaystyle{ \quad C_2:=
\supetage{J\geq 0}{n\geq N_0(J)}\|R_n^J(t_n)\|_3 \quad < \quad \infty \quad ,}
\end{array}
\end{equation} where $R_n^J = w_n^J + r_n^J$ and $N_0(J)$ is as in
Lemma \ref{clm:orth} with $\epsilon = \epsilon (J)$ chosen
appropriately.  Indeed, $C_2$ is bounded by $C_1$ and $\sup_n A_n$ (where we recall that~$\displaystyle A_n =
\sup_{0\leq t < T^*(u_{0,n})}\|u_{n}(t)\|_3$ and $A_n \to A_c <\infty$).  To show
${C_1<\infty}$, for $J$ large and $J_1 < J$ we can write
$$\|\sum_{j=0}^J \Lambda_{j,n}U_j(t_n)\|_3 \leq \sum_{j=0}^{J_1} \|\Lambda_{j,n}U_j(t_n)\|_3 + \|\sum_{j=J_1+1}^J \Lambda_{j,n}U_j(t_n)\|_3\ .$$
For $J_1$ sufficiently large, (\ref{orth3}), Lemma \ref{clm:orth} and
(\ref{ldenergy}) give a uniform bound of the second term.  Since we
stay strictly away from the potential blow-up times of all profiles,
for fixed $J_1$ the first term is bounded due to Theorem \ref{decay}
and the time-continuity in ${L^3}$ of each $U_j$.

 Set $v_n: =
u_n - \Lambda_{0,n}U_0$.  Due to (\ref{ineq}), we have
$$
\longformule{
\left| \|u_n(t_n)\|_3^3 - \| \Lambda_{0,n}U_0(t_n)\|_3^3 - \|
v_n(t_n)\|_3^3 \right| \lesssim \int
|\Lambda_{0,n}U_0(t_n)|^{2}|v_n(t_n)| }{+ \int
|\Lambda_{0,n}U_0(t_n)||v_n(t_n)|^{2}\ .}
$$
Unlike in \cite{kk} where there is only a cross-term similar to the
first one (hence one may leave the absolute value outside the integral
and use weak convergence), both terms require the use of specific
information about the components of $v_n$. We deal with the second term
first and then briefly indicate how the first one can be dealt with in
a similar way.

 We would therefore like to show that
$$\int  |\Lambda_{0,n}U_0(t_n)||v_n(t_n)|^{2} \to 0$$
as $n\to\infty$ along some subsequence.  Fix some $L_0 \in \mathbb{N} $
large, to be chosen precisely later.  Then for $J>L_0$ we write
$$v_n(t_n) = \sum_{j=1}^{L_0} \Lambda_{j,n}U_j(t_n) + \sum_{j=L_0 +1}^{J} \Lambda_{j,n}U_j(t_n) + R_n^J(t_n) =: A_1 + A_2 + A_3$$
so that
$$(v_n(t_n))^{2} = (A_1)^{2} + (A_2)^{2} + (A_3)^{2} + 2 (A_1 A_2+A_2A_3+A_1A_3)\ .$$ For the first term, note
that using the arguments in the proof of Lemma~\ref{clm:orth} we can
make $\displaystyle \int |\Lambda_{0,n}U_0(t_n)||A_1|^{2}$ arbitrarily small for
sufficiently large $n$ (depending on $L_0$, which we shall fix in a
moment) by orthogonality of the scales/cores.  For the second term,
using H\"older's inequality and Lemma \ref{clm:orth}, for an
appropriate subsequence of $n$ depending on $J$, we have
$$
\int  |\Lambda_{0,n}U_0(t_n)||A_2|^{2} \leq \|U_0\|_3
\left(2\sum_{j=L_0 +1}^{J}\left\| \Lambda_{j,n}U_j(t_n)\right\|_3^3
\right)^\frac{2}{3}
$$
which can be made arbitrarily small by choosing $L_0$ sufficiently
large due to (\ref{orth3}) and arguments similar to the proof of
(\ref{zzdprime}).  Using arguments similar to those above and noting
that $\|A_1\|_3$ and $\|A_3\|_3$ are uniformly bounded by
(\ref{consts}), we can treat all remaining terms except for those of
the form
$$
\int  |\Lambda_{0,n}U_0(t_n)||w_n^J(t_n)|^2 \quad \textrm{or} \quad
\int  |\Lambda_{0,n}U_0(t_n)||r_n^J(t_n)|^2 \,  .
$$
Since $U_0(s)\in {L^3}$, using H\"older's inequality, (\ref{orth4})
and heat estimates we can control the term involving $w_n^J$ by a
uniform constant times the quantity
$$\|\Lambda_{0,n}U_0(t_n)w_{n}^J(t_{n})\|_{\frac{3}{2}}\ .$$
Then recalling that $t_n=\lambda _{0,n}^2s$ and $s_p = -1 +
\frac{3}{p}$, approximating~$U_0(s)$ in $L^3$ by a smooth compactly
supported function and noting that one may replace $K(\tau)$ by $e^{\tau \D}$ in (\ref{besequiv}) to obtain yet another equivalent Besov norm, we can control this term by
\begin{eqnarray*}
 \|\Lambda_{0,n}U_0(t_n)\|_{p'} \|w_{n}^J(t_{n})\|_{p} &=& \|U_0(s)\|_{p'} \lambda _{0,n}^{-s_p} \|e^{t_n\Delta }\psi_{n}^J\|_{p} \\
 & \lesssim & s^{-\frac{s_p}{2}}\|U_0(s)\|_{p'}  \|\psi_{n}^J\|_{\dot B^{s_p}_{p,\infty}}
 \end{eqnarray*}
where $\frac{2}{3} = \frac{1}{p'} + \frac{1}{p}$.  As $s$ is fixed,
this term can therefore be made small for large $J$ and then $n$ due
to (\ref{orth2a}) and the continuous embedding~${\dot B^{s_p}_{p,p}}
\hookrightarrow B^{s_p}_{p,\infty}$.

We now just need to show
that $\displaystyle \int |\Lambda_{0,n}U_0(t_n)||r_n^J|^2$ can be made
arbitrarily small for large $J$ and $n$.  By a change of variables, we have
$$\int  |\Lambda_{0,n}U_0(x,t_n)||r_n^J(x,t_n)|^2\ dx = \int \left|
  U_0\left(y, s\right) \right|
\left| \Lambda^{-1}_{0,n} r_n^J(y,s) \right|^2\ dy\ .
$$
Set $\tilde r_n^J := \Lambda^{-1}_{0,n} r_n^J(s)$. Note that $\sup_{n,J}\left\|\phantom{.} |\tilde
r_n^J|^2 \right\|_{L^\frac{3}{2}} < \infty$ by (\ref{consts}),(\ref{orth4}) and linear heat estimates, and we may assume that $U_0 \in
\mathcal{C}_0^\infty$ by approximation in $L^3$.  It therefore suffices to show that~$| \tilde r_{n(J)}^J(y)|^2$
tends strongly to zero (for some increasing $n(J)$) as $J\to\infty$ in
some Banach space $B \hookrightarrow \mathscr{S}'$ and hence tends weakly to zero in
$L^\frac{3}{2}$, making this term small for large $J$ and $n=n(J)$.

 In order to do this, we claim that, since $r_n^J$
satisfies an equation of the form (\ref{remaindereq}), there exists
$N(J)\in \mathbb{N} $ defined for all $J\geq 0$ such that
\begin{equation}\label{remreg} \sup_{\scriptsize \begin{array}{c}
J\geq 0\\ n\geq N(J)
\end{array} } \|\tilde r_n^J\|_{\dot
B^{1}_{{3/2},\infty}} < \infty
\ .
\end{equation}
Let us postpone the proof of (\ref{remreg}) for a
moment, and use it to complete the proof of Claim \ref{clm:d}.

Recalling standard
product estimates in Besov spaces, we have
$$\|fg\|_{\dot B^\sigma_{{3/2},\infty}} \lesssim \|f\|_{\dot B^{s_p}_{p,\infty}} \|g\|_{\dot B^{\sigma+1}_{{3/2},\infty}} $$
for any $\sigma >0$.  Note that this is a valid application of the
product laws since $s_p + (\sigma + 1) = \frac{3}{p} + \sigma>0$ and
$s_p < \frac{3}{p}$.  Therefore we have
$$
\|({ \tilde{r}_n^J})^2\|_{\dot B^{0}_{{3/2},\infty}} \lesssim
\|{ \tilde{r}_n^J}\|_{\dot B^{s_p}_{p,\infty}}
 \|{ \tilde{r}_n^J}\|_{\dot B^{1}_{{3/2},\infty}}\lesssim
\|{\tilde{r}_n^J}\|_{\dot B^{s_p}_{p,p}}
 \|{ \tilde{r}_n^J}\|_{\dot B^{1}_{{3/2},\infty}}\ .
$$
Hence
(\ref{remainders}) along with (\ref{remreg}) imply that
$$\left\||\tilde r_{n(J)}^J|^2\right\|_{\dot B^{0}_{{3/2},\infty}} \to 0 \qquad \textrm{as} \quad J \to \infty$$
for some $n(J)$ increasing in $J$, which concludes the proof.

We now briefly return to the proof of (\ref{remreg}), which is nothing but a simple consequence of estimates on the Duhamel term in
\cite{MarcoMoi}.  Indeed, the proof of Proposition 1 in \cite{MarcoMoi} gives the estimate
\begin{equation}\label{canplan}
\|B(f,g)(t)\|_{\dot B^1_{{3/2},\infty}} \lesssim \sup_{0<\tau<t}\|fg(\tau)\|_{L^{{3/2}}}
\end{equation}
where $B(f,g) = (\partial_t - \D)^{-1}\P \nabla \cdot (f\otimes g)$ with $B(f,g)(0)=0$.
According to \eqref{remaindereq} we can write~$r_n^J$ as a sum of a
finite number (independent of $J$) of terms (each of which do depend
on $J$) of the form $B(f,g)$ which can all be controlled, after
applying \eqref{canplan} and then H\" older in $x$, by \eqref{consts},
\eqref{orth4} and standard heat estimates, plus a sum of the
form $$\quad \sum_{0\leq j \neq j' \leq J} B\left(\Lambda_{j,n}U_j , \Lambda_{j',n}U_{j'}\right)\
.$$ After applying \eqref{canplan}, we can bound this term by a
constant independent of $J$ by the orthogonality of the scales/cores
(as in the proof of Lemma \ref{clm:orth}) for $n\geq N(J)$
sufficiently large for any $J$.  Applying~$\Lambda^{-1}_{0,n}$ (under
which all norms concerned are invariant) everywhere
establishes~(\ref{remreg}) and we are done with the quadratic term
$|v_n|^2$.

We now go back to proving
$$\int  |\Lambda_{0,n}U_0(t_n)|^2|v_n(t_n)| \to 0\ ,$$
using the same decomposition of $v_n$ as a sum of three
terms. Applying the triangle inequality, terms with $A_1$ and $A_2$ go
to zero by the same arguments of orthogonality of
scales/cores. Similarly, the term in $A_3$ involving $w_n^J$ goes to
zero using the heat decay estimates. Hence all we are left with is
$$ \int \left|
  U_0\left(y, s\right) \right|^2
\left| \Lambda^{-1}_{0,n} r_n^J(y,s) \right|\ dy\ .
$$
We just proved that $r_n^J\in \dot B^{1}_{3/2,\infty}$, while we know
that $r_n^J$ goes to zero in, say, $\dot B^{-1/4}_{4,4}$, where we take $p=4$ for concreteness (general $p$ is treated similarly). By interpolation,
we get that $r_n^J\in \dot B^{1/2}_{2,20/3}$ and goes to zero in this
later norm. We conclude using composition rules in Besov spaces $\dot
B^{s}_{p,q}$, with $0<s<1$, as $\|\, |f| \, \|_{\dot B^{s}_{p,q}} \lesssim \|f\|_{\dot B^{s}_{p,q}}$ for such $s$ (a fact which readily follows from the
characterization of Besov spaces in terms of finite differences in that
range, and the elementary inequality $||a|-|b||\leq |a-b|$). As
$|U_0|^2$ is smooth, hence in the dual space $\dot
B^{-1/2}_{2,20/17}$, this last remaining integral goes to zero, and Claim \ref{clm:d} is proved. \hfill
$\Box$

\subsection{Proof of Theorem \ref{thm:ac}}
Theorem \ref{thm:ac} is a
consequence of the following lemma which is proved in the last section of \cite{kk}, following
the argument in \cite{ess}:
\begin{lemma}\label{lemma:cont} Suppose $u_0 \in {L^3}$ and $NS(u_0)
$ belongs to~$ L^\infty([0,T];{L^3})$ for some finite $T>0$.  Then there exists
some $R_0 >>1$ such that~$u$ belongs to~$ \mathcal{C}^\infty((\mathbb{R}^3
\backslash B_{R_0}(0) \times [0,T])$, with global bounds on
derivatives.
\end{lemma}
   Indeed, assuming $T^*(u_0) < +\infty$ and applying Lemma
\ref{lemma:cont} with $T=T^*(u_0)$, $u(t) \rightharpoonup 0$ as $t\nearrow T^*(u_0)$ implies
that $D^\alpha u(x,T^*) \equiv 0$ for $|x|>R_0$ for any multi-index
$\alpha$.  Then known backwards uniqueness and unique continuation
results for the parabolic inequality satisfied by $\omega := \nabla
\times u$ show that $\omega \equiv 0$ on $\mathbb{R}^3 \times
[0,T^*]$, see the last section of~\cite{kk} for more
details (see also \cite{ess}). This implies
$u\equiv 0$ as well due to the divergence-free condition, and hence
$T^*(u_0) = +\infty$ by uniqueness of mild solutions, contrary to
assumption, which proves Theorem \ref{thm:ac}.  \hfill $\Box$

\section{ Minimal Blow-up Initial Data}
\label{minimal}
   In this section we consider the question of ``minimal
blow-up initial data" in various settings, of the type addressed in
\cite{sverakrusin}.

Suppose $X=X(\rd)$ is a Banach space of initial
data on which there is a norm which is invariant under the
transformations leaving the Navier-Stokes equations invariant, with
the property that there exists some small $\epsilon _0 = \epsilon
_0(X) >0$ such that $T^*(u_0) = +\infty$ whenever $\|u_0\|_X <
\epsilon_0 $, where~$T^*(u_0)$ is the maximal time of existence of
$NS(u_0)$ in the space $X$.

Then the question to be considered
is the following:
\begin{statement}\label{statement} Suppose there exists   $v_0\in
X$ such that $T^*(v_0) < \infty$, and define $\rho = \rho_X:=\inf
\{\|v_0\|_X \ | \ T^*(v_0) < +\infty\} \geq \epsilon _0 > 0$.  Then
there exists $u_0 \in X$ such that $T^*(u_0) < \infty$, and $\|u_0\|_X
= \rho$.  Moreover, up to transformations under which the
Navier-Stokes equations are invariant, the set of such $u_0$ is
compact in $Y$, for a similar space $Y$ such that $X\subseteq
Y$.
\end{statement}
   Such a statement was proved in the case $X = {\dot
H^\frac{1}{2}}(\mathbb{R}^3)$ in \cite{sverakrusin} (in fact in the
setting of weak solutions), and moreover with $Y=X$.  In the following
we show that the result is a simple consequence of the profile
decompositions, Theorem 2 in \cite{Gallagher} for $X={\dot
H^\frac{1}{2}}$ (and in fact this can easily be extended to $X=\dot
H^{\frac{d}{2}-1}(\mathbb{R}^d)$ for any $d$), Theorems
\ref{thm:dataprofa} and \ref{thm:dataprofb} stated above (proved in
\cite{gk}) and Theorem \ref{thm:profa} stated and proved above in the
settings $X= {L^d(\mathbb{R}^d)}$ and $X={\dot
B^{s_{a}}_{a,b}}(\mathbb{R}^d)$.  To be precise, what we prove is the
following:
\begin{thm}\label{mindata}
Statement \ref{statement} is true for
$X=Y=\dot H^{\frac{d}{2}-1}(\mathbb{R}^d)$ for any~$d \geq 2$, and there
exists a norm on $X$, equivalent to the standard norm, such that
Statement \ref{statement} is true for
$(X,Y)=({L^d(\mathbb{R}^d)},{\dot B^{s_{p}}_{p,q}}(\mathbb{R}^d))$
whenever~$3\leq d < p\leq q \leq \infty$, and for $(X,Y) = ({\dot
B^{s_{a}}_{a,b}}(\mathbb{R}^d),{\dot B^{s_{p}}_{p,q}}(\mathbb{R}^d))$
for any $d\geq 3 $ and $a, b\in [1,2d+3)$ satisfying $\max\{a,b\} < p
\leq \infty$ and~$1\leq b < (p/a)b  \leq q \leq \infty$, where
$s_r:= -1 + \frac{d}{r}$ for $r\in \R$.
\end{thm}

   Note that it is important that $b<\infty$ so that local solutions
are in fact known to exist (and hence a maximal time of existence
makes sense), as opposed to only having global solutions for small
data.  (We shall see below that this is necessary for a different
technical reason as well.)  Note also that in applying Theorem \ref{thm:profa} in the proof below, one may have to use a set of smaller indices first in the space $Y$ (to satisfy the assumptions of that theorem), and then the more general statement follows from the standard embeddings.

\medskip

{\bf Proof of Theorem \ref{mindata}.} \quad For simplicity, we first prove the theorem for $X= {\dot
B^{s_{a}}_{a,b}}$, and define~$\| \cdot \|_X:= \| \cdot {\tilde
\|}_{\dot B^{s_{a}}_{a,b}}$ (this norm is defined via wavelet bases,
see~\cite{gk}). Assume there is some ${\dot B^{s_{a}}_{a,b}}$ datum
with a finite maximal time of existence, so that $\rho = \rho_{\dot
B^{s_{a}}_{a,b}}$ is well-defined.  By known small data regularity
results there exists $\epsilon _0$ such that $\rho \geq \epsilon _0
>0$.  By the definition of $\rho$, there exists a sequence $u_{0,n}
\in {\dot B^{s_{a}}_{a,b}}$ with $T^*(u_{0,n})<\infty$ (hence
necessarily $\|u_{0,n}{\tilde \|}_{\dot B^{s_{a}}_{a,b}} \geq \rho$)
and $\|u_{0,n}{\tilde \|}_{\dot B^{s_{a}}_{a,b}} \searrow \rho$ as
$n\to\infty$.

Since $u_{0,n}$ is therefore a bounded sequence in
${\dot B^{s_{a}}_{a,b}}$, we can apply the profile decomposition
Theorem \ref{thm:dataprofb} to this sequence, so that, after passing
to a subsequence, we may write $u_{0,n}$ as
$$u_{0,n}(x) = \sum_{j=0}^J \frac 1{\lambda_{j,n}}  \phi_ j\left(\frac{ x - x_{j,n}}{ \lambda_{j,n}} \right) + \psi_n ^J(x)\ ,$$
and (\ref{orth3a}) gives
\begin{equation}\label{dataorth} \sum_{j=0}^\infty \| \phi_j {\tilde
\|}_{\dot B^{s_{a}}_{a,b}}^{\tau} \leq \liminf_{n'\to\infty}
\|u_{0,n'}{\tilde \|}_{\dot B^{s_{a}}_{a,b}}^{\tau} = \rho^\tau\ .
\end{equation} Moreover, applying Theorem \ref{thm:profa}, we see that
there is at least one $j_0\in \mathbb{N}$ such that
$T^*(\phi_{j_0}) < +\infty$.  Indeed, otherwise one could take~$\tau_n
\equiv +\infty$ in that theorem and see that~$NS(u_{0,n})$ lives past
its finite maximal time of existence which is impossible.  By
definition of~$\rho$, we know that~$\|\phi_{j_0}{\tilde \|}_{\dot
B^{s_{a}}_{a,b}} \geq \rho$, else we would have~$T^*(\phi_{j_0}) =
+\infty$.  However,~(\ref{dataorth}) gives~$\|\phi_{j_0}{\tilde
\|}_{\dot B^{s_{a}}_{a,b}} \leq \rho$, so that~$\|\phi_{j_0}{\tilde
\|}_{\dot B^{s_{a}}_{a,b}} = \rho$ and we may take~$u_0 = \phi_{j_0}$
in the statement of the Theorem.  This proves the existence statement.

 \begin{remark}\label{onlyoneprofilecompact}
 We remark similarly to Remark~\ref{onlyoneprofile} in the previous section that  (\ref{dataorth}) implies that only one
profile appears in the decomposition of~$u_{0,n}$.
 \end{remark}

To prove the compactness statement, suppose now moreover that
$$\|u_{0,n}{\tilde \|}_{\dot B^{s_{a}}_{a,b}} \equiv \rho$$ (having
established the existence of at least one such element) and passing to
a subsequence write $u_{0,n}$ in a profile decomposition as before.
The same results hold as above, and note that as pointed out in Remark~\ref{onlyoneprofilecompact},   (\ref{dataorth}) implies
that~$\phi_{j} = 0$ for all~$j\neq j_0$, since necessarily~$\|\phi_{j_0}{\tilde \|}_{\dot B^{s_{a}}_{a,b}} = \rho$ implies that~$
\sum_{j\neq j_0} \| \phi_j {\tilde \|}_{\dot B^{s_{a}}_{a,b}}^{\tau}
\leq 0$. Note that here we have used the fact that~$b<\infty$ so
that $\tau < \infty$.  Therefore we can write
$$u_{0,n}(x) = \frac 1{\lambda_n}  \phi \left(\frac{ x - x_n}{ \lambda_n} \right) + \psi_n (x) =:\Lambda_n \phi(x) + \psi_n(x)\ ,$$
where $\phi = \phi_{j_0}$, etc., and $\psi_n \to 0$ in ${\dot
B^{s_{p}}_{p,q}}$ as $n\to\infty$ by (\ref{orth2a}).  The invariances
of the ${\dot B^{s_{p}}_{p,q}}$ norm imply as well that $\Lambda_n^{-1}
\psi_n \to 0$ in~${\dot B^{s_{p}}_{p,q}}$,
hence clearly $\Lambda_n^{-1} u_{0,n}  \to \phi$
in ${\dot B^{s_{p}}_{p,q}}$, and the theorem is proved.

 To prove
the theorem for $X={L^d(\mathbb{R}^d)}$, we consider a minimizing
sequence in ${L^d(\mathbb{R}^d)}$ and proceed as above applying
Theorem \ref{thm:profa}.  We similarly conclude that there exists some
profile $\phi_{j_0}$ with finite maximal time of existence in ${\dot
B^{s_{a}}_{a,b}}$.  Otherwise, $NS(u_{0,n})$ would be globally defined
in ${\dot B^{s_{a}}_{a,b}}$, and standard ``persistency" results for
Navier-Stokes (see, e.g., \cite{gip3}) would then imply that $NS(u_{0,n})$ is global in
${L^d(\mathbb{R}^d)}$ as well, contrary to assumption.  The remainder
of the proof follows as above due to Remark \ref{rkaaa} and the theorem is proved in this case as well.

For $X=\dot H^{\frac{d}{2}-1}(\mathbb{R}^d)$, the proof is identical
using the theorems in \cite{Gallagher} (with the usual norm on $X$)
and we would initially take $Y = {L^d(\mathbb{R}^d)}$ to see that $\Lambda_n^{-1}u_{0,n} \to \phi$ in $Y$.  Moreover,
since $X \hookrightarrow Y$, $\phi \in X $ and $\|\Lambda_n^{-1}u_{0,n}\|_X=\|u_{0,n}\|_X \equiv
\|\phi\|_X = \rho$, we see that $\Lambda_n^{-1}u_{0,n} \rightharpoonup \phi$ in $X$, and since $X$ is a Hilbert space the above properties imply that actually we have strong convergence in $X$ as desired.  \hfill $\Box$

\appendix
\section{A perturbation result}
\newcommand\normtb[5]{\|#5\|_{\mathcal{L}^{#1}_T\dot{B}^{#2}_{#3,#4}}}
\newcommand\normb[4]{\|#4\|_{\dot{B}^{#1}_{#2,#3}}}
\newcommand\vbar{\overline{v}}

\def\virgp{\raise 2pt\hbox{,}}

Let us  state (without proof) a perturbation result for the
$d$-dimensional Navier-Stokes system.
 \begin{prop}\label{pdr1}
  Let~$\displaystyle
s_p=-1+\frac dp
$, $r\in [1,\frac {2p}{p+1}]$ and define $\displaystyle
 s:=s_p+\frac2r$. Assume finally that~$p<2d+3$. There are constants~$\varepsilon_0$ and~$C$ such that the following holds.
Let~$w_0 \in \dot B^{s_p}_{p,p}$, $f \in F:={\mathcal L}^r([0,T];\dot
B^{s-2}_{p,p})+{\mathcal L}^{\frac{2p}{p+1}}([0,T];\dot B^{s_p-1+\frac
  1 p}_{p,q})$ and~$v \in  D:=\mathcal L^p([0,T];\dot B^{s_p+\frac2p}_{p,p})$ be given, with
$$
\normb {s_p}pp{w_0} +\|{f}\|_{F} \leq \varepsilon_0 \exp
\left( - C  \|v\|_{D} \right).
$$
Suppose moreover that $\mbox{div} \: v=0$, and let $w$ be a solution of
$$
  \partial_t w-\Delta w+w\cdot\nabla w
+v \cdot\nabla w +w \cdot\nabla v =-\nabla \pi + f
$$
with~$
\mbox{div} \: w=0.
$
Then~$w$ belongs to~$E_{p,p}(T)$ and the following
estimate holds:
$$
\|w\|_{E_{p,p}(T)}\lesssim (\normb {s_p}pq{w_0} + \|{f}\|_F) \exp C  \|{v}\|_D.
$$

\end{prop}
The proof of that proposition follows the estimates of~\cite{gip3}
(see in particular Propositions 4.1 and Theorem 3.1
of~\cite{gip3}). The two main differences are
\begin{itemize}
\item the absence of an exterior
force in~\cite{gip3}, but that force is added with no difficulty to
the estimates ;
\item the rather weak estimate on the drift term $v$, which accounts
  for the restricted numerology on time exponents in the definition of
  $E_{p,p}$. The reader should note that closing estimates on $w$ in
  our setting amounts to doing again the same estimates that were done
  in the proof of Lemma \ref{lemmaremainder}.
\end{itemize}

\end{document}